\documentclass[11pt]{amsart}

\usepackage{amsmath,amssymb}

\def\({\left(}  \def\){\right)} \def\[{\left[}  \def\]{\right]}
 \newcommand{\eproof}{\hfill $\Box$} \def\OL{\overline}
   \def
\RR{{\rm {\mathbb R}}} \def \NN{{\rm {\mathbb N}}}   
\newtheorem{theorem}{Theorem} \newtheorem{lemma}{Lemma}
\newtheorem{remark}{Remark}
\newtheorem {corollary}{Corollary}

\newtheorem{proposition}{Proposition}
\newtheorem{definition}{Definition}

\newcommand{\beq}{\begin{equation}}
\newcommand{\eeq}{\end{equation}}

\begin {document}
\begin {center}
A TOPOLOGICAL SEPARATION CONDITION FOR FRACTAL ATTRACTORS

\bigskip

T. Bedford$^a$, S. Borodachov$^b$, J. Geronimo$^{c,}$\footnote {This
author was partially supported by NSF grant DMS-0500641.}

\bigskip

{\it $^a$University of Strathclyde, Glasgow, Scotland, G1 1QX}

{\it $^b$Towson University, Baltimore, MD, 21252-0001

$^c$Georgia Institute of Technology, Atlanta, GA, 30332-0160}
\end {center}

\bigskip

{\bf Abstract.} We consider finite systems of contractive
homeomorphisms of a complete metric space, which are non-redundant
on every level. In general this separation condition is weaker than
the strong open set condition and is not equivalent to the weak
separation property. We prove that this separation condition is equivalent
to the strong Markov property (see definition below). We also show that the set of $N$-tuples of
contractive homeomorphisms, which are non-redundant on every level,
is a $G_\delta$ set in the topology of pointwise convergence of
every component mapping with an additional requirement that the
supremum of contraction coefficients of mappings be strictly less
than one. We give several sufficient conditions for this separation
property. For every fixed $N$-tuple of $d\times d$ invertible
contraction matrices from a certain class, we obtain density results
for $N$-tuples of fixed points which define $N$-tuples of mappings
non-redundant on every level.



\section {Notation and definitions.}

Let $X$ be a complete metric space and $d$ be the distance in $X$.
Recall that a mapping $w:X\to X$ is called a {\it contracting
mapping} (or a {\it contraction}) if
$$
\sigma=\sigma (w)=\sup\limits _{x\neq y\in X}{\frac {{d}(w(x),w(y))}{{d}(x,y)}}<1.
$$
The number $\sigma(w)$ will be referred to as the contraction
coefficient of the mapping $w$.

Let $N\in \NN$, $w_1,\ldots,w_N:X\to X$ be contracting
homeomorphisms of $X$ onto itself and $A=A(w_1,\ldots,w_N)\subset X$
be the unique non-empty compact set such that
$$
A=\bigcup\limits_{i=1}^{N}{w_i(A)}.
$$
The set $A$ is known as the invariant set or the attractor of the
system $\{w_1,\ldots,w_N\}$. This way to define the attractor first
appears in the paper by Hutchinson \cite {Hut81}. Denote
$\Sigma=\{1,\ldots,N\}$ and for every vector ${\bf
i}=\{i_1,\ldots,i_n\}\in \Sigma^n$, let
$$w_{{\bf
i}}=w_{i_1,\ldots,i_n}=w_{i_1} \ldots w_{i_n}=w_{i_1}\circ \ldots
\circ w_{i_n}.$$ Denote by $\mathcal M(X)$ the space of all
contracting homeomorphisms $w:X\to X$ of the space $X$ onto itself.
\begin {definition}
{\rm For every $n\in \NN$, denote by $\mathcal V_n$ the set of all
ordered $N$-tuples $(w_1,\ldots,w_N)\in (\mathcal M(X))^N$ such that
for every ${\bf i}\in \Sigma^n$, there holds
$$
w_{\bf i}(A)\nsubseteq \bigcup _{{\bf j}\in \Sigma^n,\ {\bf j}\neq
{\bf i}}{w_{\bf j}(A)}.
$$
}
\end {definition}

We say that a system $(w_1,\ldots,w_N)\in (\mathcal M(X))^N$ satisfies the open set condition (OSC), if
there is a non-empty open set $\mathcal O\subset X$ such that

1. $w_{i}(\mathcal O)\cap w_{j}(\mathcal O)=\emptyset$, $i\neq j$;

2. $w_i(\mathcal O)\subset \mathcal O$, $i=1,\ldots,N$.

We say that the system $(w_1,\ldots,w_N)$ satisfies the strong open
set condition (SOSC) if it satisfies the OSC with $\mathcal O\cap
A\neq \emptyset$.

A mapping $w:X\to X$ is called a {\it contracting similitude} if there is a number $\sigma\in (0,1)$ such that
$$
d(w(x),w(y))=\sigma d(x,y),\ \ x,y\in X.
$$
The attractor of a finite system of contracting similitudes in $X$
is known as {\it self-similar set}. When $X=\RR^d$, $d\in \NN$, and
$w_1,\ldots,w_N:\RR^d\to\RR^d$ are contracting similitudes, the SOSC
and the OSC are equivalent (cf. the result of Schief \cite {Sch94}).
In general, the OSC does not imply the SOSC (cf. e.g. \cite{Sch94}).
The above definition of self-similarity is different from the
definition given for example in the book by Mattila \cite
{MatGSMES}, where additional restrictions on the size of the
overlaps are required.

We say that a collection $(w_1,\ldots,w_N)\in (\mathcal M(X))^N$ satisfies the Markov partition property (MPP)
if there exists a subset $V\subset A$ open relatively to $A$ such
that

1. $\OL V=A$;

2. $w_i(V)\cap w_j(V)=\emptyset$, $i\neq j$.

\begin {definition}
{\rm We say that a system of mappings $(w_1,\ldots,w_N)\in (\mathcal
M(X))^N$ satisfies the strong Markov property (SMP) if for every
$n\in \NN$, there is an open set $\mathcal O_n\subset X$ such that

1. $\OL {\mathcal O_n\cap A}=A$;

2. $w_{\bf i}(\mathcal O_n)\cap w_{\bf j}(\mathcal O_n)=\emptyset$,
for every ${\bf i}\neq {\bf j}\in \Sigma^n$. }
\end {definition}

It is not difficult to see that SMP implies MPP if we let
$V=\mathcal O_1\cap A$, and that SOSC implies the SMP if we set
$\mathcal O_n=\mathcal O$ for every $n\in \NN$ (see Proposition \ref
{Th1}). The SMP does not in general imply the SOSC (see Remark \ref
{R2} below). Hence, MPP is also a weaker property than SOSC. We also
remark here that SMP does not follow from the weak separation
property of Lau and Ngai (see Example 2 on p. 76 in \cite
{LauNga99}).

Denote by $\Sigma^\infty$ the set of all infinite sequences
$(i_1,i_2,\ldots)$, where $i_j\in \Sigma$, $j=1,2,\ldots$. A
sequence $(i_1,i_2,\ldots)\in \Sigma^\infty$ is called an address of
a point $x\in A$, if
$$
x\in \bigcap\limits_{n=1}^{\infty}{w_{i_1,\ldots,i_n}(A)}.
$$
This is equivalent to the fact that for some point $a\in X$,
$$
x=\lim\limits_{n\to\infty}{w_{i_1,\ldots,i_n}({a})}.
$$
It is not difficult to see that every point $x\in A$ has at least
one address and every sequence from $\Sigma^\infty$ is an address of
some point from $A$. The set $$\mathcal T=\bigcup\limits_{i\neq
j}{w_i(A)\cap w_j(A)}$$ is non-empty if and only if there are points
in $A$, which have more than one address.

An interesting question is how generic are any of the above
separation conditions in ${\mathcal M}(X)$. One of the results we
present below is to show that the SMP condition is a countable
intersection of open sets, i.e. a $G_{\delta}$ set. This result
should be contrasted with that of Falconer~\cite{Fal88} where he
considered attractors associated with affine maps and obtained a
formula for the Hausdorf dimension that was generic in the sense of
Lebesgue measure (see also results by Mattila \cite [Theorem
9.13]{MatGSMES} and Solomyak \cite {Sol98}).

In Section \ref {s2} we show that SMP holds if and only if
$(w_1,\ldots,w_N)$ is non-redundant on every level, i.e. $(w_1,\ldots,w_N)\in\cap_{n=1}^{\infty}{\mathcal V}_n$. Furthermore
we show that the set of all systems of mappings that satisfy SMP is
a $G_\delta$ set in a suitable topology. The proofs of these results
are presented in Sections \ref {s3} and \ref {s4}. In Section \ref
{s5} we find certain sufficient conditions for the SMP. In Section
\ref {s6}, we discuss the relation between the SMP for a
self-similar set in $\RR^d$ and the equality of its similarity and
Hausdorff dimension. Section \ref {s7} deals with density results
for the SMP in the case of self-affine sets in $\RR^d$.

\section {Main results}\label {s2}

\begin {theorem}\label {equivalence}
Let $X$ be a complete metric space. The system $(w_1,\ldots,w_N)$ of
contracting homeomorphisms of $X$ onto $X$ satisfies the SMP if and
only if
$$
(w_1,\ldots,w_N)\in \bigcap\limits_{n=1}^{\infty}{\mathcal V_n}.
$$
\end {theorem}

\begin {definition}
{\rm We will call a sequence $\{w^m\}_{m\in \NN}$ from $\mathcal
M(X)$ strongly pointwise convergent to a mapping $w\in \mathcal
M(X)$ and write $w^m{{\ s.p. \over}\!\!\!\!\!\to } w$, $m\to\infty$,
if

1. $\lim\limits_{m\to\infty}{w^m(x)}=w(x)$ for every $x\in X$;

2. $\sup\limits_{m\in \NN}{\sigma (w^m)}<1$.}
\end {definition}

If $\{w^m\}_{m\in \NN}\subset \mathcal M(X)$ is a sequence of
similitudes and $w\in M(X)$ is a similitude, then strong pointwise
convergence is equivalent to the ``usual" pointwise convergence.

We introduce a topology $\mathcal B_N$ on the space $(\mathcal
M(X))^N$ by defining a subset $C\subset (\mathcal M(X))^N$ to be
closed if for every sequence
$\{(w^m_1,\ldots,w^m_N)\}_{m\in\NN}\subset C$, such that
$\{w^m_i\}{{\ s.p. \over}\!\!\!\!\!\to }w_i\in \mathcal M(X)$,
$i=1,\ldots,N$, we have $(w_1,\ldots,w_N)\in C$. We agree here that
$\emptyset$ is closed. It is not difficult to see, for example, that
the space $(\mathcal M(X))^N$ with the topology $\mathcal B_N$ is a
Hausdorff space.

\begin {theorem}\label {G-delta}
Let $N\in \NN$ and $X$ be a complete metric space. The set of systems of mappings $(w_1,\ldots,w_N)\in (\mathcal M(X))^N$, which satisfy the SMP is a $G$-delta set in the topology $\mathcal B_N$.
\end {theorem}

For a $d\times k$ matrix $B$, let
\begin {equation}\label {norm}
\|B\|=\max\limits_{{\bf x}\in \RR^k\setminus \{{\bf 0}\}}{\frac
{\left|B{\bf x}\right|}{\left|{\bf x}\right|}}
\end {equation}
be its norm. We say that $B$ is a {\it contraction matrix} if
$\|B\|<1$.

Let $X=\RR^d$ and $B_1,\ldots,B_N$ be invertible $d\times d$
contraction matrices. Denote by $E_d(B_1,\ldots,B_N)$ the set of all
ordered $N$-tuples
$(\boldsymbol\alpha_1,\ldots,\boldsymbol\alpha_N)$ of points from
$\RR^d$ such that the system of mappings $w_i:\RR^d\to \RR^d$,
$$
w_i({\bf x})=B_i({\bf x}-\boldsymbol\alpha_i)+\boldsymbol\alpha_i,\
\ \ i=1,\ldots,N,
$$
satisfies the SMP. We will sometimes consider the set
$E_d(B_1,\ldots,B_N)$ as a subset of $\RR^{dN}$.
\begin {corollary}\label {C3}
For any collection $B_1,\ldots,B_N$ of invertible $d\times d$
contraction matrices, the set $E_d(B_1,\ldots,B_N)$ is a $G$-delta
subset of $\RR^{dN}$ (in the topology induced by the Euclidean
distance).
\end {corollary}

\section {Proof of Theorem \ref {equivalence}}\label {s3}

We will start the proof with the following statement.
\begin {lemma}\label {MPP}
Let $X$ be a complete metric space and $(w_1,\ldots,w_N)\in (\mathcal M(X))^N$. If $(w_1,\ldots,w_N)\in \cap_{n=1}^{\infty}{\mathcal V_n}, $ then
there is an open set $\mathcal O\subset X$ such that $\OL {\mathcal
O\cap A}=A$ and $w_i(\mathcal O)\cap w_j(\mathcal O)=\emptyset$,
$i\neq j$. In particular, the system $(w_1,\ldots,w_N)$ will satisfy
the MPP.
\end {lemma}

{\bf Proof.} In order to prove Lemma \ref {MPP} denote
$$K_i(A)=w_i(A)\setminus \bigcup\limits _{j=1\atop j\neq
i}^{N}{w_j(A)},\ \ i=1,\ldots,N.
$$
Let also $$Z_i=w_i^{-1}(K_i(A))\ \ {\rm and}\ \
V=\bigcap\limits_{i=1}^{N}{Z_i}.$$ For example, if $w_1(x)=x/2$ and
$w_2(x)=x/2+1/2$, then $A=[0,1]$, $Z_1=[0,1)$, $Z_2=(0,1]$, and
hence, $V=(0,1)$.

It is not difficult to see that $Z_i\subset A$, $i=1,\ldots,N$. We
show that $\OL Z_i=A$, $i=1,\ldots,N$. Let $x\in A$ and let
$U\subset X$ be any open set containing $x$. Denote by $B(a,\rho)$
the open ball in $X$ centered at point $a$ of radius $\rho>0$. Since
$w_i(U)$ is also open, there is $\epsilon
>0$ such that $B(w_i(x),\epsilon)\subset w_i(U)$. Let $r_i=\sigma(w_i)\in (0,1)$
be the contraction coefficient of $w_i$, $i=1,\ldots,N$, and define
$$r_{\rm max}=\max\limits_{i=1,\ldots,N}{r_i}.$$
Choose a number $m\in \NN$ so that $(r_{\rm max})^m\cdot {\rm
diam}A<\epsilon$. There exist indices $i_1,\ldots,i_m\in \Sigma$
such that $x\in w_{i_1,\ldots,i_m}(A)$. Then $w_i(x)\in
w_{i,i_1,\ldots,i_m}(A)$ and
$$
{\rm diam}\ w_{i,i_1,\ldots,i_m}(A)\leq r_i\cdot r_{i_1}\cdot \ldots
\cdot r_{i_m}\cdot{\rm diam}A\leq (r_{\rm max})^{m+1}\cdot{\rm
diam}A<\epsilon.
$$
Hence,
\begin {equation}\label {1}
w_{i,i_1,\ldots,i_m}(A)\subset B(w_i(x),\epsilon)\subset w_i(U).
\end {equation}
Since $(w_1,\ldots,w_N)\in \mathcal V_{m+1}$, we have
$$
w_{i,i_1,\ldots,i_m}(A)\nsubseteq
\bigcup\limits_{j_1,\ldots,j_{m+1}\in \Sigma\atop j_1\neq
i}{w_{j_1,\ldots,j_{m+1}}(A)}=\bigcup\limits_{j=1\atop j\neq
i}^{N}{w_j(A)}.
$$
Hence, there is $z\in A$ such that $w_{i,i_1,\ldots,i_m}(z)$ does
not belong to $\cup_{j: j\neq i}\ {w_j(A)}.$ Let
$t=w_{i_1,\ldots,i_m}(z)$. Since $w_{i}(t)$ does not belong to any
$w_j(A)$ with $j\neq i$, we must have $w_{i}(t)\in w_i(A)$, that is,
$w_i(t)\in K_i(A)$. Hence, $t\in Z_i$. On the other hand, since
$w_i(t)\in w_{i,i_1,\ldots,i_m}(A)$, in view of (\ref {1}), we have
$w_i(t)\in w_i(U)$, that is $t\in U$, which implies that $\OL
Z_i=A$, $i=1,\ldots,N$.

We next show that $\OL V=A$. Indeed, since each $Z_i$ is open
relative to $A$, there are open sets $W_i\subset X$ such that
$Z_i=W_i\cap A$, $i=1,\ldots,N$. Let $y$ be any element in $A$ and
$U$ be any open neighborhood of $y$. Since $\OL Z_1=A$, there is
$z_1\in Z_1\cap U=A\cap W_1\cap U$. Since $\OL Z_2=A$, there is
$z_2\in Z_2$ in the open neighborhood $W_1\cap U$ of the point
$z_1\in A$, that is $z_2\in A\cap U\cap W_1\cap W_2$. Then by
induction, there will be an element $z_N\in A\cap U\cap W_1\cap
\ldots\cap W_N=V\cap U$, and the required relation follows.

Note that for every $i\neq j$, there holds
$$
w_i(V)\cap w_j(V)\subset w_i(Z_i)\cap w_j(Z_j)=K_i(A)\cap
K_j(A)\subset
$$
$$
\subset\(w_i(A)\setminus w_j(A)\)\cap w_j(A)=\emptyset.
$$
Taking also into account the fact that $V$ is relatively open with
respect to $A$ as an intersection of a finite collection of subsets
of $A$, which are open relative to $A$, we conclude that the system
$(w_1,\ldots,w_N)$ possesses the MPP.

For every $x\in V$, denote
$$
\rho (x)=\min_{i=1,\ldots,N}{{\rm dist}\(w_i(x),\bigcup
\limits_{j=1\atop j\neq i}^{N}{w_j(A)}\)}.
$$
In view of the relations
$$
w_i(V)\subset w_i(Z_i)=K_i(A),\ \ i=1,\ldots,N,
$$
point $w_i(x)$, $x\in V$, does not belong to the closed set $\bigcup
\limits_{j: j\neq i}\ {w_j(A)}$. Hence, $\rho (x)>0$, $x\in V$, and
the set
$$
\mathcal O=\bigcup \limits_{x\in V}{B\(x,\rho (x)/2\)}
$$
is open. Since $\OL V=A$ and $V\subset \mathcal O\cap A\subset A$,
we have $\OL {\mathcal O\cap A}=A$. To show that $w_i(\mathcal
O)\cap w_j(\mathcal O)=\emptyset$, $i\neq j$, assume to the contrary
that there exist indices $i\neq j$ such that $w_i(\mathcal O)\cap
w_j(\mathcal O)$ contains some element $y$. Then $y=w_i(p)=w_j(q)$
for some $p,q\in \mathcal O$. There are points $c,b\in V$ such that
$d(c,p)<\rho(c)/2$ and $d(b,q)<\rho(b)/2$. Note that
\begin {equation}\label {3}
d(y,w_i(c))=d(w_i(p),w_i(c))\leq r_i\cdot d(p,c)<r_i\cdot \rho(c)/2
\end{equation}
and
\begin {equation}\label {2}
d(y,w_j(b))=d(w_j(q),w_j(b))\leq r_j\cdot d(q,b)<r_j\cdot \rho(b)/2.
\end {equation}
There also hold the following relations
\begin {equation}\label {4}
\rho(c)\leq {\rm dist}\(w_i(c),\bigcup\limits_{k=1\atop k\neq
i}^{N}{w_k(A)}\)\leq {\rm dist}(w_i(c),w_j(A))\leq d(w_i(c),w_j(b))
\end {equation}
and
\begin {equation}\label {5}
\rho(b)\leq {\rm dist}\(w_j(b),\bigcup\limits_{k=1\atop k\neq
j}^{N}{w_k(A)}\)\leq {\rm dist}(w_j(b),w_i(A))\leq d(w_j(b),w_i(c)).
\end {equation}
Then, in view of relations (\ref {3})--(\ref {5}), we obtain
$$
\rho(c)+\rho(b)\leq 2d(w_i(c),w_j(b))\leq
2(d(w_i(c),y)+d(y,w_j(b)))<
$$
$$
<r_i\cdot \rho(c)+r_j\cdot \rho (b)<\rho (c)+\rho (b),
$$
which is impossible. Hence, $w_i(\mathcal O)$ and $w_j(\mathcal O)$
are disjoint, which completes the proof of Lemma~\ref{MPP}.\eproof

To prove sufficiency in Theorem \ref {equivalence}, assume that
$$
(w_1,\ldots,w_N)\in \bigcap\limits_{n=1}^{\infty}{\mathcal
V_n}\subset (\mathcal M(X))^N.
$$
Then for every $m\in \NN$ and $n\in \NN$, we have $(w_1,\ldots,w_N)\in \mathcal V_{nm}\subset (\mathcal M(X))^N$, which implies that the system $\{w_{\bf i}\}_{{\bf i}\in
\Sigma^m}$ belongs to the set
${\mathcal V_n}\subset (\mathcal M(X))^{N^m}$. Hence, $\{w_{\bf i}\}_{{\bf i}\in
\Sigma^m}\in \cap_{n=1}^{\infty}{\mathcal V_n}\subset (\mathcal M(X))^{N^m}$. By Lemma \ref {MPP}, there is an open set $\mathcal
O_m\subset X$ such that $\OL {\mathcal O_m \cap A}=A$ and $w_{\bf
i}(\mathcal O_m)\cap w_{\bf j}(\mathcal O_m)=\emptyset$ for every
${\bf i}\neq {\bf j}\in \Sigma^m$, $m\in \NN$. Hence, the system
$(w_1,\ldots,w_N)$ satisfies the SMP.

The proof of the necessity in Theorem \ref {equivalence} is preceded by the following
proposition.
\begin {lemma}\label {open}
Let mappings $w_1,\ldots,w_N\in \mathcal M(X)$ be such that there is
a non-empty open set $\mathcal O\subset X$ with the property
$$
w_i({\mathcal O})\cap w_{j}(\mathcal O)=\emptyset, \ \ \ i\neq j.
$$
Then for every $i=1,\ldots,N$,
$$
w_{i}(\mathcal O)\cap \bigcup \limits_{j : j\neq i}{w_j(\OL
{\mathcal O})}=\emptyset.
$$
\end {lemma}
{\bf Proof.} Assume the contrary. Then for some $j_0\neq i$, there
$x\in w_i(\mathcal O)\cap w_{j_0}(\OL {\mathcal O})$. Let $z\in \OL
{\mathcal O}$ be such that $x=w_{j_0}(z)$. There is a sequence
$\{z_m\}_{m\in \NN}\subset \mathcal O$ such that
$z=\lim\limits_{m\to\infty}{z_m}$ and hence,
$x=\lim\limits_{m\to\infty}{w_{j_0}(z_m)}$. Since $w_i(\mathcal O)$
is an open neighborhood of $x$, we have $w_{j_0}(z_m)\in
w_i(\mathcal O)$ for every $m$ sufficiently large, and hence,
$w_i(\mathcal O)\cap w_{j_0}(\mathcal O)\neq \emptyset$, which
contradicts the assumption, thus Lemma \ref {open} is proved.\eproof

{\bf Completion of the proof of Theorem \ref{equivalence}.} Assume
that system $(w_1,\ldots,w_N)\in (\mathcal M(X))^N$ satisfies the
SMP. Let $k\in \NN$ be arbitrary. Then there is an open set
$\mathcal O_k\subset X$ such that $\OL {\mathcal O_k \cap A}=A$ and
$w_{\bf i}(\mathcal O_k)\cap w_{\bf j}(\mathcal O_k)=\emptyset$ for
every ${\bf i}\neq {\bf j}\in \Sigma^k$. We show that for every
${\bf i}\in \Sigma^k$,
\begin {equation}\label {p5}
w_{\bf i}(A)=\OL {w_{\bf i}(\mathcal O_k)\cap A}.
\end {equation}
Taking into account Lemma \ref {open} and the fact that $A=\OL
{\mathcal O_k\cap A}\subset \OL {\mathcal O_k}$, we obtain
$$
w_{\bf i}(\mathcal O_k)\cap A=\(w_{\bf i}(\mathcal O_k)\cap w_{\bf
i}(A)\)\cup \(w_{\bf i}(\mathcal O_k)\cap \bigcup\limits_{{\bf j}\in
\Sigma^k,\ {\bf j}\neq {\bf i}}{w_{\bf j}(A)}\)\subset
$$
$$
\subset w_{\bf i}(\mathcal O_k\cap A)\cup \(w_{\bf i}(\mathcal
O_k)\cap \bigcup\limits_{{\bf j}\in \Sigma^k,\ {\bf j}\neq {\bf
i}}{w_{\bf j}(\OL {\mathcal O_k})}\)=w_{\bf i}(\mathcal O_k\cap A).
$$
Then
$$
\OL {w_{\bf i}(\mathcal O_k)\cap A}\subset \OL {w_{\bf i}(\mathcal
O_k\cap A)}={w_{\bf i}(\OL {\mathcal O_k\cap A})}=w_{\bf i}(A).
$$
On the other hand,
$$
w_{\bf i}(A)={w_{\bf i}(\OL {\mathcal O_k\cap A})}=\OL {w_{\bf
i}(\mathcal O_k\cap A)}=\OL {w_{\bf i}(\mathcal O_k)\cap w_{\bf
i}(A)}\subset \OL {w_{\bf i}(\mathcal O_k)\cap A},
$$
and (\ref {p5}) follows.

Assume that $(w_1,\ldots,w_N)$ does not belong to
$\cap_{n=1}^{\infty}\mathcal V_n$. Then there is $n\in \NN$ and
${\bf i}_n\in \Sigma^n$ such that
$$
w_{{\bf i}_n}(A)\subset \bigcup\limits_{{\bf j}\in \Sigma^n,\ {\bf
j}\neq {\bf i}_n}{w_{\bf j}(A)}.
$$
Then, taking into account (\ref {p5}) we obtain
$$
w_{{\bf i}_n}(\mathcal O_n)\cap A\subset \OL {w_{{\bf i}_n}(\mathcal
O_n)\cap A}=w_{{\bf i}_n}(A)\subset \bigcup\limits_{{\bf j}\in
\Sigma^n,\ {\bf j}\neq {\bf i}_n}{w_{\bf j}(A)}=
$$
$$
=\bigcup\limits_{{\bf j}\in \Sigma^n,\ {\bf j}\neq {\bf i}_n}{\OL
{w_{\bf j}(\mathcal O_n)\cap A}}\subset \bigcup\limits_{{\bf j}\in
\Sigma^n,\ {\bf j}\neq {\bf i}_n}{w_{\bf j}(\OL{\mathcal O_n})}.
$$
Since $\OL {w_{{\bf i}_n}(\mathcal O_n)\cap A}=w_{{\bf i}_n}(A)\neq
\emptyset$, there is a point $x\in w_{{\bf i}_n}(\mathcal O_n)\cap
A\subset w_{{\bf i}_n}(\mathcal O_n)$. Then $x\in
\bigcup\limits_{{\bf j}\in \Sigma^n,\ {\bf j}\neq {\bf i}_n}{w_{\bf
j}(\OL {\mathcal O_n})}$, Hence, $$w_{{\bf i}_n}(\mathcal O_n)\cap
\bigcup\limits_{{\bf j}\in \Sigma^n,\ {\bf j}\neq {\bf i}_n}{w_{\bf
j}(\OL {\mathcal O_n})}\neq \emptyset,$$ which contradicts to Lemma
\ref {open}. Theorem \ref {equivalence} is proved.\eproof

\section {Proof of Theorem \ref {G-delta}.}\label {s4}

The proof of some statements in this section is standard, but we
include it for the convenience of the reader.
\begin {lemma}\label {1'}
If a sequence $\{w^m\}_{m\in \NN}\subset \mathcal M(X)$ converges strongly pointwise to a mapping $w\in \mathcal M(X)$, then the sequence of fixed points of mappings $w^m$ converges to the fixed point of $w$.
\end {lemma}
{\bf Proof.} Let $x^m\in X$ be the fixed point of $w^m$, $m\in \NN$, and $x\in X$ be the fixed point of $w$. Denote also
$$
\sigma =\sup\limits_{m\in \NN}{\sigma (w^m)}.
$$
Then
$$
d(x^m,x)\leq d(x^m,w^m(x))+d(w^m(x),x)=
$$
$$
=d(w^m(x^m), w^m(x))+d(w^m(x),w(x))\leq \sigma d(x^m,x)+d(w^m(x),w(x)).
$$
Hence,
$$
d(x^m,x)\leq \frac {1}{1-\sigma}d(w^m(x),w(x)),
$$
and we have
$$
\lim\limits_{m\to\infty}{d(x^m,x)}=0.
$$
Lemma \ref {1'} is proved.\eproof

\begin {lemma}\label {2'}
Let $A$ be the attractor of a system of mappings $w_1,\ldots,w_N\in \mathcal M(X)$ with contraction coefficients not exceeding a given number $\sigma \in (0,1)$. Let also $B[a,r]$ be a closed ball containing the fixed point of every mapping $w_1,\ldots,w_N$. Then $A\subset B[a,R]$, where $R=\frac {1+\sigma}{1-\sigma} r$.
\end {lemma}
{\bf Proof.} Assume the contrary. Denote by $y_1,\ldots,y_N$ the
fixed points of mappings $w_1,\ldots,w_N$ respectively. Let $z$ be a
point in $A$ furthest from $a$. Then we must have $d(z,a)>R$. Let
$1\leq i\leq N$ be such index that $z=w_i(z_1)$ for some $z_1\in A$.
Then
$$
d(z_1,a)\geq d(z_1,y_i)-d(y_i,a)\geq \frac {1}{\sigma}d(w_i(z_1),w_i(y_i))-r=
$$
$$
=\frac {1}{\sigma}d(z,y_i)-r\geq \frac {1}{\sigma}d(z,a)-\frac {1}{\sigma}d(y_i,a)-r\geq \frac {1}{\sigma}d(z,a)-\frac {r}{\sigma}-r.
$$
Hence,
$$
\frac {d(z_1,a)}{d(z,a)}\geq\frac {1}{\sigma}-\frac {(1+\sigma)r}{\sigma d(z,a)}>\frac{1}{\sigma}-\frac {(1+\sigma)r}{\sigma R}=1,
$$
which contradicts to the fact that $z$ is a point in $A$ furthest
from $a$.\eproof

\begin {lemma}\label {3'}
Let $\{w^m_1\}_{m\in \NN},\ldots, \{w^m_n\}_{m\in \NN}$ be sequences of mappings from $\mathcal M(X)$ such that $w^m_i{{\  s.p. \over}\!\!\!\!\!\to \  }w_i\in \mathcal M(X)$, $i=1,\ldots,n$. Then $w^m_1\circ \ldots \circ w^m_n {{\  s.p. \over}\!\!\!\!\!\to\ }w_1\circ \ldots \circ w_n$, $m\to \infty$.
\end{lemma}
{\bf Proof.} We will use induction. For $n=1$, the assertion of the
lemma is trivial. Assume that the assertion is true for a given
value of $n\geq 1$ and show that it holds for any $n+1$ sequences
satisfying the assumptions of the lemma. For every $x\in X$, we will
have
\begin{align*}
d(w^m_1 w^m_2 \ldots  w^m_{n+1}(x),w_1 w_2  \ldots w_{n+1}(x))
&\leq d(w^m_1(w^m_2 \ldots  w^m_{n+1}(x)),w^m_1(w_2  \ldots  w_{n+1}(x)))\\
&\quad +
d(w^m_1(w_2 \ldots w_{n+1}(x)),w_1(w_2  \ldots  w_{n+1}(x)))\\
&\leq d(w^m_2 \ldots  w^m_{n+1}(x),w_2 \ldots w_{n+1}(x))\\
&\quad +d(w^m_1(w_2 \ldots w_{n+1}(x)),w_1(w_2 \ldots  w_{n+1}(x))).
\end{align*}
By the assumption of the induction, both distances in the last line
vanish as $m\to\infty$ and we have
$$
\lim\limits_{m\to\infty}{w^m_1 w^m_2 \ldots w^m_{n+1}(x)}=w_1 w_2 \ldots  w_{n+1}(x),\ \ x\in X.
$$
Since $$\sigma=\max_{i=1,\ldots,n+1}{\sup\limits_{m\in \NN}{\sigma (w^m_i)}}<1,$$ we have
$$
\sigma (w^m_1 w^m_2 \ldots w^m_{n+1})\leq \sigma ^{n+1}<1,\ \ m\in \NN,
$$
which implies strong pointwise convergence. Lemma \ref {3'} is
proved.\eproof

Given a system $W=(w_1,\ldots,w_N)\in (\mathcal M(X))^N$ and an address ${\bf i}\in \Sigma^\infty$, let $\Pi_{\bf i}(W)$ be the point in the attractor of $W$ with address ${\bf i}$.
\begin {lemma}\label {4'}
Let $W_m=\(w^m_1,\ldots,w^m_N\)$, $m\in \NN$, be a sequence from $(\mathcal M(X))^N$ such that for every $i=1,\ldots,N$, the sequence $\{w^m_i\}_{m\in \NN}$ converges strongly pointwise to some mapping $w_i\in \mathcal M(X)$. Then for every address ${\bf i}\in \Sigma^\infty$,
$$
\lim\limits_{m\to\infty}{\Pi_{\bf i}(W_m)}=\Pi_{\bf i}(W),
$$
where $W=(w_1,\ldots,w_N)$.
\end {lemma}
{\bf Proof.} Given an arbitrary address ${\bf i}=(i_1,i_2,\ldots)\in \Sigma^\infty$, denote by $x_{i_1\ldots i_n}$ the fixed point of the mapping $w_{i_1\ldots i_n}$. Let also
$$
\delta=\max\limits_{i=1,\ldots,N}{\sup\limits_{m\in \NN}{\sigma
(w^m_i)}}.
$$
Let $B(a,r)$ be a ball containing the attractor $A$ of the system
$W$ and $R=\frac {1+\delta}{1-\delta}r$.

Choose an arbitrary $\epsilon >0$ and let $n\in \NN$ be large enough
so that
\begin {equation}\label {p}
d(\Pi_{\bf i}(W),x_{i_1\ldots i_n})<\epsilon\ \ \ {\rm and}\ \ \ R\delta^n <\epsilon.
\end {equation}
Denote by $x^m_{\alpha_1\ldots \alpha_n}$ the fixed point of the
mapping $w^m_{\alpha_1}\circ \ldots \circ w^m_{\alpha_n}$,
$\alpha_1,\ldots,\alpha_n\in \Sigma$. By Lemma \ref {3'}, we have
$w^m_{\alpha_1}\circ \ldots \circ w^m_{\alpha_n}{{\  s.p.
\over}\!\!\!\!\!\to\  }w_{\alpha_1\ldots \alpha_n}$, $m\to\infty$.
Then by Lemma \ref {1'}, we have
$\lim\limits_{m\to\infty}{x^m_{\alpha_1\ldots
\alpha_n}}=x_{\alpha_1\ldots \alpha_n}$, for every
$\alpha_1,\ldots,\alpha_n\in \Sigma$. Since $x_{\alpha_1\ldots
\alpha_n}\in A\subset B(a,r)$, there is a number $m_n\in \NN$ such
that for every $m>m_n$ and $\alpha_1,\ldots,\alpha_n\in \Sigma$, we
have $x^m_{\alpha_1\ldots \alpha_n}\in B(a,r)$. For every $m>m_n$,
we obtain
\begin{align*}
d(\Pi_{\bf i}(W),\Pi_{\bf i}(W_m))&\leq d(\Pi_{\bf i}(W),x_{i_1\ldots i_n})+d(w_{i_1\ldots i_n }(x_{i_1\ldots i_n}),w^m_{i_1}\ldots w^m_{i_n}(x_{i_1\ldots i_n}))\\
&\quad +d(w^m_{i_1}\ldots w^m_{i_n}(x_{i_1\ldots i_n}),\Pi_{\bf i}(W_m))\\
&\leq \epsilon +d(w_{i_1\ldots i_n }(x_{i_1\ldots i_n}),w^m_{i_1}\ldots w^m_{i_n}(x_{i_1\ldots i_n}))\\
&\quad +d(w^m_{i_1}\ldots w^m_{i_n}(x_{i_1\ldots i_n}),w^m_{i_1}\ldots w^m_{i_n}(z_{{\bf i},m})),
\end{align*}
where $z_{{\bf i},m}$ is some point in the attractor $A_m$ of the system $W_m$. Taking into account Lemma \ref {3'}, we will have
$$
d(\Pi_{\bf i}(W),\Pi_{\bf i}(W_m))\leq \epsilon +o(1)+ \delta ^n d(x_{i_1\ldots i_n},z_{{\bf i},m}).
$$
For every $i=1,\ldots,N$, the fixed point $x^m_i$ of $w^m_i$ is also the fixed point of the $n$-th power of $w^m_i$, and as it was noted above, $x^m_i\in B(a,r)$, $m>m_n$. By Lemma \ref {2'}, we have $z_{{\bf i},m}\in A_m\subset B[a,R]$.
Since $x_{i_1\ldots i_n}\in A\subset B(a,r)\subset B[a,R]$, in view of (\ref {p}), we obtain
$$
d(\Pi_{\bf i}(W),\Pi_{\bf i}(W_m))\leq \epsilon +o(1)+2R\delta^n\leq
3\epsilon +o(1).
$$
Hence,
$$
\limsup\limits_{m\to\infty}{d(\Pi_{\bf i}(W),\Pi_{\bf i}(W_m))}\leq 3\epsilon.
$$
In view of arbitrariness of $\epsilon$, we have
$$
\lim\limits_{m\to\infty}{d(\Pi_{\bf i}(W),\Pi_{\bf i}(W_m))}=0,
$$
and the assertion of Lemma \ref {4'} follows.\eproof

Let
$$
\mathcal F=\bigcup\limits_{n\in \NN}{\Sigma^n}.
$$
\begin {lemma}\label {5'}
Let $N,n\in \NN$. If a sequence $\{(w^m_1,\ldots,w^m_N)\}_{m\in
\NN}\subset (\mathcal M(X))^N\setminus \mathcal V_n$ converges
strongly pointwise in every component to a system
$W=(w_1,\ldots,w_N)\in (\mathcal M(X))^N$, then we have $W\in
(\mathcal M(X))^N\setminus \mathcal V_n$.
\end {lemma}
From Lemma \ref {5'} we obtain the following statement, which in view of Theorem \ref {equivalence}, implies the assertion of Theorem \ref {G-delta}.
\begin {corollary}\label {7}
For every positive integers $n$ and $N$, the set $ \mathcal V_n$ is open in the topology $\mathcal B_N$, and hence, $\cap_{n=1}^{\infty}{\mathcal V_n}$ is a G-delta set.
\end {corollary}
{\bf Proof of Lemma \ref {5'}.} Let $W_m=(w^m_1,\ldots,w^m_N)\in (\mathcal M(X))^N\setminus \mathcal V_n$ be a sequence, where every component is convergent strongly pointwise to the corresponding component of the system $W=(w_1,\ldots,w_N)\in (\mathcal M(X))^N$. Denote $w^m_{\bf k}=w^m_{k_1}\circ \ldots \circ w^m_{k_p}$, ${\bf k}=(k_1,\ldots,k_p)\in \mathcal F$. For every $m\in \NN$, there is a vector ${\bf i}_m\in \Sigma^n$ such that
$$
w^m_{{\bf i}_m}(A_m)\subset \bigcup\limits_{{\bf k}\in \Sigma^n,\ {\bf k}\neq {\bf i}_m}{w^m_{\bf k}(A_m)},
$$
where $A_m$ is the attractor of the system $W_m$. There is an index ${\bf i}\in \Sigma^n$ and infinite subsequence $\mathcal N\subset \NN$ such that
\begin {equation}\label {p2}
w^m_{{\bf i}}(A_m)\subset \bigcup\limits_{{\bf k}\in \Sigma^n,\ {\bf k}\neq {\bf i}}{w^m_{\bf k}(A_m)},\ \ m\in \mathcal N.
\end {equation}
Let $A$ be the attractor of the system $W$ and $x\in w_{\bf i}(A)$
be an arbitrary point. Then $x=\Pi_{{\bf i}\beta}(W)$ for some
$\beta\in \Sigma^\infty$. In view of (\ref {p2}), for every $m\in
\mathcal N$, there holds
$$
\Pi_{{\bf i}\beta}(W_m)\in w^m_{\bf i}(A_m)\cap w^m_{{\bf j}_m}(A_m)
$$
for some ${\bf j}_m\in \Sigma^n$ distinct from ${\bf i}$. There are index ${\bf j}\in \Sigma^n$, ${\bf j}\neq {\bf i}$, and infinite subsequence $\mathcal N'\subset \mathcal N$ such that
$$
\Pi_{{\bf i}\beta}(W_m)\in w^m_{\bf i}(A_m)\cap w^m_{{\bf j}}(A_m),\ \ \ m\in \mathcal N'.
$$
Hence, there is a sequence $\gamma_m=(\gamma^m_1,\gamma^m_2,\ldots)\in \Sigma^\infty$ such that
\begin {equation}\label {p3}
\Pi_{{\bf i}\beta}(W
_m)=\Pi_{{\bf j}\gamma_m}(W_m), \ \ \ m\in \mathcal N'.
\end {equation}
One can find an infinite subsequence $\mathcal N_1\subset \mathcal
N'$ and an index $\gamma_1\in \Sigma$ such that
$\gamma^m_1=\gamma_1$, $m\in \mathcal N_1$. One can find an infinite
subsequence $\mathcal N_2\subset \mathcal N_1$ and an index
$\gamma_2\in \Sigma$ such that $\gamma^m_1=\gamma_1$ and
$\gamma^m_2=\gamma_2$, $m\in \mathcal N_2$. Continuing this process
indefinitely, we obtain an address
$\gamma=(\gamma_1,\gamma_2,\ldots)\in \Sigma^\infty$ and a sequence
of embedded infinite sets $\mathcal N_1\supset \mathcal N_2\supset
\ldots \supset \mathcal N_k \supset \ldots$ such that
$\gamma^m_k=\gamma_k$, $m\in \mathcal N_k$, $k\in \NN$.

Let as above, $B(a,r)$ be an open ball containing $A$. Since $A$ contains the fixed points of mappings $w_1,\ldots,w_N$, by Lemma \ref {1'}, there is $m_0\in \NN$ such that for every $m\in \NN$, $m>m_0$, the fixed points of mappings $w^m_1,\ldots,w^m_N$ will be in $B(a,r)$. Then, by Lemma \ref {2'}, $A_m\subset B[a,R]$, where $R=\frac {1+\delta}{1-\delta}r$ and $\delta=\max\limits_{i=1,N}{\sup\limits_{m\in \NN}{\sigma (w^m_i)}}\in (0,1).$ For every $k\in \NN$ and $m\in \mathcal N_k$, $m>m_0$, there are points $b$ and $c$ in $A_m$ such that
\begin{align*}
d(\Pi_{{\bf j}\gamma_m}(W_m),\Pi_{{\bf j}\gamma}(W_m))&=d(w^m_{\bf j}w^m_{\gamma^m_1\ldots \gamma^m_k}(b),w^m_{\bf j}w^m_{\gamma_1\ldots \gamma_k}(c))\\
&=d(w^m_{\bf j}w^m_{\gamma_1\ldots \gamma_k}(b),w^m_{\bf
j}w^m_{\gamma_1\ldots \gamma_k}(c)) \leq \sigma (w^m_{\bf j})\cdot
\sigma (w^m_{\gamma_1})\cdot \ldots \cdot \sigma
(w^m_{\gamma_k})d(b,c)\\
&\leq\delta ^{n+k}{\rm diam}A_m\leq 2R\delta^{n+k}.
\end{align*}
By Lemma \ref {4'} and relation (\ref {p3}), for every $m\in \mathcal N_k$, $m>m_0$, we obtain
\begin{align*}
d(x,\Pi_{{\bf j}\gamma}(W))&\leq d(\Pi_{{\bf i}\beta}(W),\Pi_{{\bf i}\beta}(W_m))+d(\Pi_{{\bf j}\gamma_m}(W_m),\Pi_{{\bf j}\gamma}(W_m))\\
&\quad +d(\Pi_{{\bf j}\gamma}(W_m),\Pi_{{\bf j}\gamma}(W))\leq 2R\delta ^{n+k}+o(1).
\end{align*}
Hence, letting $m\to\infty$ along the sequence $\mathcal N_k$, we will have
$$
d(x,\Pi_{{\bf j}\gamma}(W))\leq 2R\delta^{n+k}, \ \ k\in \NN.
$$
Letting now $k\to\infty$ we get that $d(x,\Pi_{{\bf j}\gamma}(W))=0$, which implies that
$$
x=\Pi_{{\bf j}\gamma}(W)\in w_{\bf j}(A)\subset \bigcup\limits_{{\bf k}\in \Sigma^n,\ {\bf k}\neq {\bf i}}{w_{\bf k}(A)},
$$
where vector ${\bf j}$ was chosen to be distinct from {\bf i}.
Since $x\in w_{\bf i}(A)$ was chosen arbitrarily, we obtain that
$$
w_{\bf i}(A)\subset \bigcup\limits_{{\bf k}\in \Sigma^n,\ {\bf k}\neq {\bf i}}{w_{\bf k}(A)},
$$
and hence, $W\in (\mathcal M(X))^N\setminus \mathcal V_n$. Lemma
\ref {5'} is proved, which completes the proof of Theorem \ref
{G-delta}.\eproof

{\bf Proof of Corollary \ref {C3}.} Let $U_n$, $n\in \NN$, be the
set of ordered $N$-tuples
$(\boldsymbol\alpha_1,\ldots,\boldsymbol\alpha_N)\in (\RR^d)^N$ such
that the system of mappings
\begin {equation}\label{x'y'}
w_i({\bf x})=B_i({\bf x}-\boldsymbol\alpha_i)+\boldsymbol\alpha_i,\
\ i=1,\ldots,N,
\end{equation}
belongs to $\mathcal V_n$. By Theorem \ref {equivalence}, we have
$$E_d(B_1,\ldots,B_N)=\bigcap\limits_{n=1}^{\infty}{U_n}.$$

It remains to show that for every $n\in \NN$, the set $U_n$ is open.
Assume the contrary and let ${\boldsymbol
\alpha}=(\boldsymbol\alpha_1,\ldots,\boldsymbol\alpha_N)\in U_n$ be
not an interior point of $U_n$. Then there is a sequence
$\{{\boldsymbol \beta}_m\}_{m=1}^{\infty}\subset
(\RR^{d})^N\setminus U_n$ such that ${\boldsymbol
\alpha}=\lim\limits_{m\to\infty}{{\boldsymbol\beta}_m}$. Let
${\boldsymbol\beta}_m=(\boldsymbol\beta^m_1,\ldots,\boldsymbol\beta^m_N)$,
$m\in \NN$, where $\boldsymbol\beta^m_i\in \RR^d$, $i=1,\ldots,N$.
Then for every $m\in \NN$, the system of contracting mappings
$$
w^m_i({\bf x})=B_i({\bf
x}-\boldsymbol\beta^m_i)+\boldsymbol\beta^m_i,\ \ i=1,\ldots,N,
$$
does not belong to $\mathcal V_n$. Since for every $i=1,\ldots,N$
and ${\bf x}\in \RR^d$,
$$
\lim\limits_{m\to\infty}{w^m_i({\bf x})}=w_i({\bf x}),
$$
where $w_i$ is defined as in (\ref {x'y'}), and
$$
\max\limits_{i=1,\ldots,N}{{\|B_i\|}}<1,
$$
we have a strong pointwise convergence of the sequence
$\{w^m_i\}_{m=1}^{\infty}$ to $w_i$, $i=1,\ldots,N$. By Lemma \ref
{5'}, we have that $(w_1,\ldots,w_N)$ does not belong to $\mathcal
V_n$, i.e. $\boldsymbol\alpha\notin U_n$. This contradiction shows
that $U_n$ is an open set for every $n$ and the assertion of
Corollary \ref {C3} follows.\eproof

\section {Sufficient conditions for the SMP}\label {s5}

\begin {proposition}\label {Th1}
Let $X$ be a complete metric space and $(w_1,\ldots,w_N)$ be a
collection of contracting homeomorphisms of $X$ onto $X$. If
$(w_1,\ldots,w_N)$ satisfies the SOSC, then
$$
(w_1,\ldots,w_N)\in \bigcap\limits_{n=1}^{\infty}{\mathcal V_n},
$$
or, equivalently, $(w_1,\ldots,w_N)$ satisfies the SMP.
\end {proposition}
The converse is not true (see Remark \ref {R2} below).

{\bf Proof.} Let $\mathcal O\subset X$ be the open set from the definition of the SOSC. Show that $\OL {A\cap \mathcal O}=A$. Indeed, if $x\in A$ and $\epsilon >0$ are arbitrary, for some $m\in \NN$ sufficiently large and ${\bf i}\in \Sigma^m$ we have $x\in w_{\bf i}(A)\subset B(x,\epsilon)$. Since $w_{\bf i}(\mathcal O\cap A)\neq \emptyset$, $w_{\bf i}(\mathcal O\cap A)\subset B(x,\epsilon)$, and
$$
w_{\bf i}(\mathcal O\cap A)=w_{\bf i}(\mathcal O)\cap w_{\bf i}(A)\subset \mathcal O\cap A,
$$
we have $(\mathcal O\cap A)\cap B(x,\epsilon)\neq \emptyset$. Hence, $A\subset\OL {\mathcal O\cap A}$. Since the opposite inclusion is trivial, we have $\OL {\mathcal O\cap A}=A$.

If for every $n\in \NN$, we let $\mathcal O_n=\mathcal O$, then condition 1 in the definition of the SMP holds. For every ${\bf i}\neq {\bf j}\in \Sigma^n$, if $1\leq k\leq n$ is the smallest index such that $i_k\neq j_k$, then
\begin{align*}
w_{\bf i}(\mathcal O_n)\cap w_{\bf j}(\mathcal O_n)&=
w_{\bf i}(\mathcal O)\cap w_{\bf j}(\mathcal O)\subset w_{i_1\ldots i_{k-1}}(w_{i_k}(\mathcal O))\cap w_{i_1\ldots i_{k-1}}(w_{j_k}(\mathcal O))\\
&=w_{i_1\ldots i_{k-1}}(w_{i_k}(\mathcal O)\cap w_{j_k}(\mathcal O))=\emptyset.
\end{align*}
Thus, the system $(w_1,\ldots,w_N)$ satisfies the SMP and, by
Theorem \ref {equivalence}, we have $(w_1,\ldots,w_N)\in
\cap_{n=1}^{\infty}{\mathcal V_n}$. Proposition \ref {Th1} is
proved.\eproof

Recall that $$\mathcal T=\bigcup\limits_{i\neq j}{w_i(A)\cap w_j(A)}$$ and denote
$$
\mathcal D=\mathcal D(w_1,\ldots,w_N)=A\setminus
\bigcup\limits_{{\bf i}\in \mathcal F}{w_{\bf i}^{-1}(\mathcal T)}.
$$
The following result holds.
\begin {proposition}\label {P2}
Let $X$ be a complete metric space and $w_1,\ldots,w_N$ be
contracting homeomorphisms of the space $X$ onto $X$. If $\mathcal
D=\mathcal D(w_1,\ldots,w_N)\neq \emptyset$, then

1. $w_i(\mathcal D)\subset \mathcal D$, $i=1,\ldots,N$;

2. $w_i(\mathcal D)\cap w_j(\mathcal D)=\emptyset$, $i\neq j$;

3. $\OL {\mathcal D}=A$;

4. $(w_1,\ldots,w_N)\in \bigcap\limits_{n=1}^{\infty}{\mathcal
V_n}$.
\end {proposition}
{\bf Proof.} To prove the first statement assume the contrary, i.e.
for some $1\leq k\leq N$, there is $y\in w_{k}(\mathcal D)\setminus
\mathcal D$. Then there is a point $x\in \mathcal D$ such that
$y=w_k(x)$. On the other hand, since $y$ is not in $\mathcal D$,
there is a vector ${\bf p}=(p_1,\ldots,p_s)\in \mathcal F$ such that
$w_{\bf p}(y)\in \mathcal T$. Hence, $w_{p_1,\ldots,p_s,k}(x)\in
\mathcal T$, which contradicts to the fact that $x\in \mathcal D$.

To prove the second statement, assume again the contrary, i.e. for
some indexes $1\leq i\neq j\leq N$, there is a point $x\in
w_i(\mathcal D)\cap w_j(\mathcal D)$. Then $x=w_i(t)$, $t\in
\mathcal D$. Since
$$
w_i(t)\in w_i(\mathcal D)\cap w_j(\mathcal D)\subset w_i(A)\cap
w_j(A)\subset \mathcal T,
$$
we have a contradiction with the fact that $t\in \mathcal D$.

To show the third statement, choose any point $z\in A$ and a ball
$B(z,\epsilon)$, $\epsilon >0$. Denote $r_{\rm
max}=\max\limits_{i=1,\ldots,N}{\sigma(w_i)}.$ Let $m\in \NN$ be
such number that $r_{\max}^m\cdot {\rm diam}A<\epsilon$ and ${\bf
i}=(i_1,\ldots,i_m)\in \Sigma^m$ be such that $z\in w_{\bf i}(A)$.
Let point $q\in A$ be such that $z=w_{\bf i}(q)$ and $x$ be some
point in $\mathcal D$. Then, by the first statement, $w_{\bf
i}(x)\in \mathcal D$. Since
$$
d(z,w_{\bf i}(x))=d(w_{\bf i}(q),w_{\bf i}(x))\leq
\sigma(w_{i_1})\cdot \ldots \cdot \sigma(w_{i_m})\cdot d(q,x)\leq
r^m_{\max}\cdot {\rm diam}A<\epsilon,
$$
we have $\mathcal D\cap B(z,\epsilon)\neq \emptyset$ for every $z\in
A$ and $\epsilon >0$. Taking into account that $\mathcal D\subset
A$, we have $\OL {\mathcal D}=A$.

Statement 4 is also proved by contradiction. Assume that $\mathcal
D\neq \emptyset$, but $(w_1,\ldots,w_N)$ does not belong to
$\mathcal V_n$ for some $n\in \NN$. Then there is a vector ${\bf
i}=(i_1,\ldots,i_n)\in \Sigma^n$ such that
$$
w_{\bf i}(A)\subset \bigcup\limits_{{\bf j}\in \Sigma^n\atop {\bf
j}\neq {\bf i}}{w_{\bf j}(A)}.
$$
Let $x$ be any point in $\mathcal D$. There is a vector ${\bf k}\in
\Sigma^n$, ${\bf k}=(k_1,\ldots,k_n)\neq {\bf i}$, such that $w_{\bf
i}(x)\in w_{\bf i}(A)\cap w_{\bf k}(A)$. If $i_1\neq k_1$, then
$w_{\bf i}(x)\in w_{i_1}(A)\cap w_{k_1}(A)\subset \mathcal T$, which
contradicts to the fact that $x\in \mathcal D$.

If $i_1=k_1$, let $1\leq s<n$ be an index that $i_1=k_1$, ...,
$i_s=k_s$, but $i_{s+1}\neq k_{s+1}$. Then
\begin{align*}
w_{\bf i}(x)&\in w_{i_1,\ldots,i_s}\(w_{i_{s+1},\ldots,i_n}(A)\)\cap
w_{i_1,\ldots,i_s}\(w_{k_{s+1},\ldots,k_n}(A)\)\\
&=w_{i_1,\ldots,i_s}\(w_{i_{s+1},\ldots,i_n}(A)\cap
w_{k_{s+1},\ldots,k_n}(A)\)\\
&\subset w_{i_1,\ldots,i_s}\(w_{i_{s+1}}(A)\cap
w_{k_{s+1}}(A)\)\subset w_{i_1,\ldots,i_s}\(\mathcal T\).
\end{align*}
Hence, $w_{i_{s+1},\ldots,i_n}(x)\in \mathcal T$, which again
implies that $x$ does not belong to $\mathcal D$. Thus, our
assumption is wrong and the fourth statement holds. Proposition \ref
{P2} is proved.\eproof

The following statement shows the relation between the cardinality of the overlaps of sets $w_i(A)$ and the SMP.
\begin {proposition}\label {countable}
Let $w_1,\ldots,w_N\in \mathcal M(X)$ be such that the corresponding attractor $A$ is uncountable and every set $w_i(A)\cap w_j(A)$, $i\neq j$, is at most countable. Then the system $(w_1,\ldots,w_N)$ satisfies the SMP.
\end {proposition}
{\bf Proof.} By assumption, the set $\mathcal T$ is at most
countable. Then the set $\cup_{{\bf i}\in \mathcal F}{w^{-1}_{\bf
i}(\mathcal T)}$ is also at most countable. Since $A$ is
uncountable, we have $D(w_1,\ldots,w_N)\neq \emptyset$. By
Proposition \ref {P2}, we have $(w_1,\ldots,w_N)\in
\cap_{n=1}^{\infty}{\mathcal V_n}$, which in view of Theorem \ref
{equivalence}, implies the SMP. Proposition \ref {countable} is
proved.\eproof

\begin {proposition}\label {P5}
Let $X$ be a complete metric space and $w_1,\ldots,w_N:X\to X$ be
contracting homeomorphisms of $X$ onto $X$. Assume that every point
in the attractor $A$ of this system has a finite number of
addresses. Then
$$
(w_1,\ldots,w_N)\in \bigcap\limits_{n=1}^{\infty}{\mathcal V_n},
$$
or equivalently, the system $(w_1,\ldots,w_N)$ satisfies the SMP.
\end {proposition}
{\bf Proof.} Assume the contrary, i.e. there exist $n\in \NN$ and
${\bf i}\in \Sigma^n$ such that
\begin {equation}\label {a1}
w_{\bf i}(A)\subset \bigcup\limits_{{\bf j}\in \Sigma^n ,\ {\bf
j}\neq {\bf i}}{w_{\bf j}(A)}.
\end {equation}
Denote by $x$ the fixed point of $w_{\bf i}$. In view of (\ref
{a1}), $x\in w_{{\bf j}}(A)$ for some ${\bf j}\in \Sigma^n$, ${\bf
j}\neq {\bf i}$. Let $p\in A$ be such point that $x=w_{{\bf j}}(p)$.
Then for every $m\in \NN$, we have $x=(w_{\bf i})^m(x)=(w_{\bf i})^m
w_{{\bf j}}(p)$. Hence, for every $m\in \NN$, the point $x$ will
have an address starting with $m$ vectors ${\bf i}$ followed by
vector ${\bf j}$ different from ${\bf i}$, thus having infinitely
many addresses, which contradicts our assumption. Proposition
\ref {P5} is proved.\eproof

We say that two vectors ${\bf i},{\bf j}\in \mathcal F$ are {\it
incomparable} if neither ${\bf i}$ is an initial word of ${\bf j}$
nor ${\bf j}$ is an initial word of ${\bf i}$. Denote
$$
\mathcal E=\{w^{-1}_{\bf j}w_{\bf i} : {\bf i},{\bf j}\in \mathcal
F,\ {\bf i},{\bf j} \ {\rm incomparable}\}.
$$
Denote by $I$ the identity mapping from $X$ to $X$. In the case when
$X=\RR^d$ and $w_i$'s are contractive similitudes, the results of
papers by Hutchinson \cite {Hut81}, Bandt and Graf \cite {BanGra92},
and Schief \cite {Sch94} imply that SOSC is equivalent to the
condition that $I\notin\OL {\mathcal E}$ in the topology of
pointwise convergence of similitudes. The weak separation property
(WSP) introduced by Lau and Ngai in \cite {LauNga99} was shown to be
equivalent to the condition that $I\notin\OL {\mathcal E\setminus
\{I\}}$ for a wide class of self-similar sets (cf. the work by
Zerner \cite {Zer96}). We see that SOSC does not allow $I\in
\mathcal E$. The WSP allows $I$ to be in $\mathcal E$ as an isolated
point. The following proposition shows the relation between the
condition that $I\notin\mathcal E$ and the SMP.
\begin {proposition}\label {P5'}
Let $X$ be a complete metric space and let the system
$(w_1,\ldots,w_N)\in (\mathcal M(X))^N$ satisfy the SMP. Then
$I\notin\mathcal E$. The converse is not true.
\end {proposition}
\begin {remark}\label {R1}
{\rm This proposition together with above mentioned results implies that
for a wide class of systems of contracting similitudes in $\RR^d$,
SMP together with WSP is equivalent to SOSC.}
\end {remark}

{\bf Proof.} Assume that $I\in \mathcal E$. Then $I=w^{-1}_{\bf
j}w_{\bf i}$ for some incomparable ${\bf i}, {\bf j}\in \mathcal F$.
Hence, $w_{\bf j}=w_{\bf i}$. Without loss of generality we can
assume that vector-index ${\bf i}$ is of the same or of a shorter
length than ${\bf j}$. Since ${\bf i}$ is not a prefix of ${\bf j}$,
we have
$$
w_{\bf i}(A)=w_{\bf j}(A)\subset \bigcup\limits_{{\bf k}\neq {\bf
i}\atop \left|{\bf k}\right|=\left|{\bf i}\right|}{w_{\bf k}(A)},
$$
which implies that the SMP does not hold. Hence, SMP implies that
$I\notin\mathcal E$.

The following counterexample shows that the converse is not true.
Let $w_1(x)=x/2$, $w_2(x)=(x+1)/2$, and $w_3(x)=(x+a)/2$, where $a$
is an irrational number from $(0,1)$. It is not difficult to see
that interval $[0,1]$ is the attractor of the system
$(w_1,w_2,w_3)$. Since $w_3([0,1])\subset w_1([0,1])\cup
w_2([0,1])$, the system $(w_1,w_2,w_3)$ does not satisfy the SMP. If
we assumed that $I\in \mathcal E$, there would be incomparable
indexes ${\bf i}=(i_1,\ldots,i_n),{\bf j}=(j_1,\ldots,j_m)\in
\mathcal F$ such that $w_{\bf i}=w_{\bf j}$. Hence, $$w_{\bf
i}(x)=\frac {1}{2^n}x+\sum\limits_{k=1\atop i_k=2}^{n}{\frac
{1}{2^k}}+\sum\limits_{k=1\atop i_k=3}^n{\frac {1}{2^k}a}=w_{\bf
j}(x)=\frac {1}{2^m}x+\sum\limits_{k=1\atop j_k=2}^m{\frac
{1}{2^k}}+\sum\limits_{k=1\atop j_k=3}^m{\frac {1}{2^k}a}.$$ Then
$n=m$ and
$$
a\(\sum\limits_{k\ \! :\ \! j_k=3}{\frac {1}{2^k}}-\sum\limits_{k\
\! :\ \! i_k=3}{\frac {1}{2^k}}\)=\sum\limits_{k\ \! :\ \!
i_k=2}{\frac {1}{2^k}}-\sum\limits_{k\ \! :\ \! j_k=2}{\frac
{1}{2^k}}.
$$
Since $a$ is irrational, we must have
$$
\sum\limits_{k\ \! :\ \! j_k=3}{\frac {1}{2^k}}=\sum\limits_{k\ \!
:\ \! i_k=3}{\frac {1}{2^k}}.
$$
Hence, $\{k : i_k=3\}=\{k : j_k=3\}$. But then
$$
\sum\limits_{k\ \! :\ \! i_k=2}{\frac {1}{2^k}}=\sum\limits_{k\ \!
:\ \! j_k=2}{\frac {1}{2^k}}.
$$
Hence, $\{k : i_k=2\}=\{k : j_k=2\}$. This implies that $\{k :
i_k=1\}=\{k : j_k=1\}$ and ${\bf i}={\bf j}$, which contradicts to
the incomparability of ${\bf i}$ and ${\bf j}$. This contradiction
shows that for the system $(w_1,w_2,w_3)$ we have $I\notin \mathcal
E$ but SMP does not hold.\eproof

\section {Some results for self-similar sets in $\RR^d$.}\label {s6}

An address $(i_1,i_2,\ldots)\in \Sigma^\infty$ is called {\it
universal}, if for any vector ${\bf j}=(j_1,\ldots,j_s)\in \mathcal
F$, there is $k\geq 0$ such that $i_{k+1}=j_1$, ..., $i_{k+s}=j_s$.
An address $(i_1,i_2,\ldots)\in \Sigma^\infty$ is called {\it
recurrent}, if for every $n\in \NN$, there is $k\in \NN$ such that
$i_{k+1}=i_1$, ..., $i_{k+n}=i_n$. In other words, a universal
address is an address, which contains every finite sequence of
numbers from $\Sigma$, and a recurrent address is an address where
any finite prefix occurs further in that address.

Recall that a mapping $w:\RR^d\to\RR^d$, $d\in \NN$, is a contracting
similitude, if there is a number $r\in (0,1)$ such that
$d(w(x),w(y))=r\cdot d(x,y)$, $x,y\in \RR^d$. Here $d$ will denote
the Euclidian distance in $\RR^d$. The attractor of a system of
finite contracting similitudes is called a self-similar set.

We will need the following result.
\begin {theorem}\label {BR}(Bandt and Rao \cite {BanRao07}).
Let $w_1,\ldots,w_N:\RR^d\to \RR^d$, $d\in \NN$, be a system of
contracting similitudes and $A$ be the attractor of this system. If
one point $a\in w_{s_1}(A)$ with a recurrent address
$(s_1,s_2,\ldots)$ belongs to the set $w_{t_1}(A)$ with $t_1\neq
s_1$, then the OSC cannot hold.
\end {theorem}
We obtain the following statement in the case of self-similar
attractors in $\RR^d$.
\begin {proposition}\label {P3}
Let $X=\RR^d$, $d\in \NN$, and mappings $w_1,\ldots,w_N$ be
contracting similitudes in $\RR^d$. If
$D(w_1,\ldots,w_N)=\emptyset$, then the system $(w_1,\ldots,w_N)$
does not satisfy the OSC.
\end {proposition}
{\bf Proof.} Let $x\in A$ be a point with a universal address ${\bf
j}=(j_1,j_2,\ldots)\in \Sigma^\infty$. Since $\mathcal D=\emptyset$,
we have
$$
x\in A\subset \bigcup \limits_{{\bf i}\in \mathcal F}{w_{\bf
i}^{-1}(\mathcal T)}.
$$
Hence, there exists a vector ${\bf i}=(i_1,\ldots,i_k)\in \mathcal
F$ such that $w_{\bf i}(x)\in \mathcal T$. This imlies that there is
an index $l\in \Sigma$, $l\neq i_1$, such that $w_{\bf i}(x)\in
w_{i_1}(A)\cap w_{l}(A)$. The address ${\bf i}{\bf
j}=(i_1,\ldots,i_k,j_1,j_2,\ldots)\in \Sigma^\infty$ is also
universal. Since
$$
w_{\bf i}(x)=w_{\bf
i}\(\lim\limits_{n\to\infty}{w_{j_1,\ldots,j_n}({\bf
0})}\)=\lim\limits_{n\to\infty}{w_{i_1,\ldots,i_k,j_1,\ldots,j_n}({\bf
0})},
$$
sequence ${\bf i}{\bf j}$ is an address of the point $w_{\bf i}(x)$.
Since every universal address is also a recurrent address, by
Theorem \ref {BR}, the OSC does not hold for the system
$(w_1,\ldots,w_N)$.\eproof

Let $W=(w_1,\ldots,w_N)$, where $w_1,\ldots,w_N:\RR^d\to \RR^d$,
$d\in \NN$, are similitudes with similarity coefficients
$r_1,\ldots,r_N\in (0,1)$ respectively. Denote by $\alpha=\alpha
(W)$ the unique positive number such that
$$
r_1^\alpha +\ldots + r_N^\alpha =1.
$$
This number is known as the similarity dimension of the attractor
$A$ associated with the system $W$. Denote by ${\rm dim}A$ the
Hausdorff dimension of the set $A$ and by $\mathcal H_\lambda$,
$\lambda>0$, the $\lambda$-dimensional Hausdorff measure in $\RR^d$.
The standard covering argument shows that
\begin {equation}\label {5''}
{\rm dim}A(W)\leq \alpha (W).
\end{equation}
\begin {proposition}\label {P4}
Let $W=(w_1,\ldots,w_N)$ be a system of contracting similitudes in
$\RR^d$, $d\in \NN$, and ${\rm dim}A(W)=\alpha (W)$. Then $A$ satisfies the SMP.
\end {proposition}
\begin {remark}\label {R2}
{\rm The results by Hutchinson \cite[Theorem 1, Section 5.3] {Hut81} combined with the
results by Schief \cite[Theorem 2.1] {Sch94}  imply that for the attractor
$A(W)$ of a finite system $W$ of contracting similitudes in $\RR^d$,
we have $ \mathcal H_{\alpha (W)}{(A(W))}>0 $ if and only if $W$
satisfies the OSC. Proposition \ref {P4} implies that any finite
system of contracting similitudes $W$ such that ${\rm
dim}A(W)=\alpha (W)$ and $\mathcal H_{\alpha (W)}(A(W))=0$
(existence of such systems is proved by Solomyak \cite{Sol98}), will
still belong to $\cap_{n=1}^{\infty}{\mathcal V_n}$ and in view of
Theorem \ref {equivalence}, will have the SMP. But such system will
not satisfy the OSC. This disproves the conjecture about the
equivalence of these two properties. Since SMP implies MPP as
asserted by Lemma \ref {MPP}, we conclude that MPP is also weaker
than OSC.}
\end{remark}

{\bf Proof of Proposition \ref {P4}.} Assume the contrary. Then in
view of Theorem \ref {equivalence}, there is $n\in \NN$ such that
$(w_1,\ldots,w_N)$ does not belong to $\mathcal V_n$. Then there is
a vector ${\bf i}\in \Sigma^n$ such that
$$
w_{\bf i}(A)\subset \bigcup\limits_{{\bf j}\in \Sigma^n,\ {\bf
j}\neq {\bf i}}{w_{\bf j}(A)}.
$$
Hence,
$$
A=\bigcup\limits_{{\bf j}\in \Sigma^n}{w_{\bf
j}(A)}=\bigcup\limits_{{\bf j}\in \Sigma^n,\ {\bf j}\neq {\bf
i}}{w_{\bf j}(A)}
$$
and $A$ will be also the attractor for the system of mappings
$S=\{w_{\bf j}\}_{{\bf j}\in \Sigma^n,\ {\bf j}\neq {\bf i}}$. In
this case the similarity dimension of $A$ associated with system $S$
satisfies
$$
\sum\limits_{{\bf j}\in \Sigma^n,\ {\bf j}\neq {\bf i}}{r_{\bf
j}^{\alpha(S)}}=1,
$$
where $r_{\bf j}$ is the contraction coefficient of the mapping
$w_{\bf j}$, ${\bf j}\in \Sigma^n$. Since
$$
\sum\limits_{{\bf j}\in \Sigma^n}{r_{\bf j}^{\alpha(W)}}=1,
$$
we have $\alpha(S)<\alpha (W)$. Then, by (\ref {5''}), we obtain
${\rm dim}A\leq \alpha (S)<\alpha (W)$, which contradicts to the
assumptions of the proposition. Proposition \ref {P4} is
proved.\eproof

\section {Density of the SMP on certain classes of self-similar
sets}\label {s7}

Let $B_1,\ldots,B_N$ be invertible $d\times d$ contraction matrices,
$d\in \NN$. Recall that $E_d(B_1,\ldots,B_N)$ is the set of ordered
point collections $\(\boldsymbol \alpha_1,\ldots,\boldsymbol
\alpha_N\)\in \(\RR^d\)^N$ such that the system of mappings
$$
u_i({\bf x})=B_i({\bf x}-\boldsymbol \alpha_i)+\boldsymbol
\alpha_i,\ \ i=1,\ldots,N,
$$
has the SMP. We will consider the set $E_d(B_1,\ldots,B_N)$ as a
subset of $\RR^{dN}$.
\begin {remark}\label {R3}
{\rm When matrices $B_1,\ldots,B_N$ are orthogonal and

1) $\|B_1\|<\frac 12$, ..., $\|B_N\|<\frac 12$,

2) $\sum\limits_{i=1}^{N}{\|B_i\|^d}<1$,\\
the set $E_d(B_1,\ldots,B_N)$ is a subset of $\RR^{dN}$ of full
measure. This follows from results of Falconer \cite [Theorem
5.3]{Fal88}, Solomyak \cite [Proposition 3.1]{Sol98}), and
Proposition \ref {P4}. (Recent results of Falconer and Miao \cite {Fal08}
imply an upper estimate for the Hausdorff dimension of the
complement of $E_d(B_1,\ldots,B_N)$.) Hence, $E_d(B_1,\ldots,B_N)$
will be dense in $\RR^{dN}$. In this paper we can show that
$E_d(B_1,\ldots,B_N)$ is dense when assumption 1) is replaced with
certain other assumptions. }
\end {remark}

\begin {theorem}\label {Th5}
Let $B_1,\ldots,B_N$ be invertible $d\times d$ contraction matrices
such that $\sum\limits_{i=1}^{N}{\|B_i\|}<1$. Then the set
$E_d(B_1,\ldots,B_N)$ is a dense $G_\delta$-subset of $\RR^{dN}$.
\end {theorem}
When $d=1$ the density result of Theorem \ref {Th5} immediately follows
from the result of Mattila \cite [Theorem 9.13]{MatGSMES}. He also 
mentions without proof 
that his result can be extended to certain cases of contractive
multidimensional similitudes.
\begin {theorem}\label {Th7}
Let $B_i=\sigma_iU_i$, where $\sigma_i\in (0,1)$, $U_i$ is a
$2\times 2$ rotation matrix, $i=1,\ldots,N$, and
$\sum_{i=1}^{N}{\sigma_i^2}<1$. Then the set $E_2(B_1,\ldots,B_N)$
is either empty or is a dense $G_\delta$-subset of $\RR^{2N}$.
\end {theorem}
\begin {remark}\label {R4}
{\rm The set $E_2(B_1,\ldots,B_N)$ can be empty under assumptions of
Theorem \ref {Th7} as the following example shows. Let
$\sigma_1,\sigma_2>0$ be such that $\sigma_1+\sigma_2>1$ and
$\sigma_1^2+\sigma_2^2<1$, and $B_1=\sigma_1I_2$, $B_2=\sigma_2I_2$
(here and below $I_d$ denotes the $d\times d$ identity matrix). For
any ordered pair $(\boldsymbol \alpha_1,\boldsymbol \alpha_2)$ of
points in $\RR^2$, the attractor $A$ of the system of mappings
$$
w_i({\bf x})=B_i({\bf x}-\boldsymbol \alpha_i)+\boldsymbol
\alpha_i=\sigma_i {\bf x}+(1-\sigma_i)\boldsymbol\alpha_i,\ \ i=1,2,
$$
is the closed segment with endpoints $\boldsymbol \alpha_1$ and
$\boldsymbol \alpha_2$. The set $w_1(A)\cap w_2(A)$ is a segment of
positive length. For $n\in \NN$ sufficiently large and some index
${\bf i}\in \Sigma^n$, there holds $w_{\bf i}(A)\subset w_1(A)\cap
w_2(A)$. If ${\bf i}$ starts with $1$, we have
$$
w_{\bf i}(A)\subset w_2(A)=\bigcup\limits_{{\bf j}\in
\Sigma^{n-1}}{w_2w_{\bf j}(A)}\subset \bigcup\limits_{{\bf j}\in
\Sigma^n\atop {\bf j}\neq {\bf i}}{w_{\bf j}(A)}.
$$
If ${\bf i}$ starts with $2$ we use analogous argument. Thus, the
system $(w_1,w_2)$ does not posses the SMP for any collection of
fixed points $(\boldsymbol \alpha_1,\boldsymbol \alpha_2)$ and
hence, $E_2(B_1,B_2)=\emptyset$.}
\end {remark}

The proof of Theorems \ref {Th5} and \ref {Th7} will follow from the
statement presented below. For an ordered collection of points
$\boldsymbol \beta =\(\boldsymbol \beta_1,\ldots,\boldsymbol
\beta_N\)\in (\RR^d)^N$, denote by $\Pi_{\bf k}(\boldsymbol\beta)$
the element with the address ${\bf k}\in \Sigma^\infty$ in the
attractor of the system of mappings
$$
u_i({\bf x})=B_i({\bf x}-\boldsymbol \beta_i)+\boldsymbol \beta_i,\
\ i=1,\ldots,N.
$$
\begin {proposition}\label {P1}
Let $1\leq k\leq d$ be integers and $B_1,\ldots,B_N$ be invertible
$d\times d$ contraction matrices such that
$\sum\limits_{i=1}^{N}{\|B_i\|^k}<1$. Assume that there is an
ordered collection $\boldsymbol \gamma_1=\(\boldsymbol
\gamma^1_1,\ldots,\boldsymbol \gamma ^1_N\)\in (\RR^d)^N$ such that
the system $W=(w_1,\ldots,w_N)$, where
$$
w_i({\bf x})=B_i({\bf x}-\boldsymbol \gamma ^1_i)+\boldsymbol
\gamma^1_i,\ \ \ i=1,\ldots,N,
$$
has the SMP. In the case $k\geq 2$ assume also that there are
collections $\boldsymbol \gamma _j=\(\boldsymbol
\gamma^j_1,\ldots,\gamma^j_N\)\in (\RR^d)^N$, $j=2,\ldots,k$, such
that for every pair of addresses ${\bf i}\neq {\bf j}\in
\Sigma^\infty$ such that $\Pi_{\bf i}(\boldsymbol \gamma_1)\neq
\Pi_{\bf j}(\boldsymbol \gamma_1)$, the system of vectors
$\{\Pi_{\bf i}(\boldsymbol \gamma_i)-\Pi_{\bf j}(\boldsymbol
\gamma_i) : i=1,\ldots,k\}$ is linearly independent.

Then the set $E_d(B_1,\ldots,B_N)$ is a dense $G_\delta$-subset of
$\RR^{dN}$.
\end {proposition}
{\bf Proof.} Let ${\boldsymbol \alpha}=\(\boldsymbol
\alpha_1,\ldots,\boldsymbol \alpha_N\)\in (\RR^d)^N$ be arbitrary.
For every ${\bf t}=(t_1,\ldots,t_k)\in \RR^k$, denote by $W_{\bf
t}=(w^{\bf t}_1,\ldots,w^{\bf t}_N)$ the system of mappings
$$
w^{\bf t}_i({\bf x})=B_i({\bf x}-\boldsymbol \alpha_i-t_1\boldsymbol
\gamma^1_i-\ldots -t_k\boldsymbol \gamma^k_i)+\boldsymbol
\alpha_i+t_1\boldsymbol \gamma^1_i+\ldots +t_k\boldsymbol
\gamma^k_i,\  i=1,\ldots, N .
$$
Let $A_{\bf t}=A(W_{\bf t})$ be the attractor of the system $W_{\bf
t}$ and $A=A(W)$ be the attractor of the system $W$. Denote
$$
P(\boldsymbol \alpha)=\{{\bf t}\in \RR^k : W_{\bf t} \ {\rm has \
no\ SMP}\}
$$
and for an index ${\bf i}=(i_1,\ldots,i_n)\in \Sigma^n$, let
$$
w^{\bf t}_{\bf i}=w^{\bf t}_{i_1}\circ \ldots \circ w^{\bf t}_{i_n}.
$$
Then
$$
P(\boldsymbol
\alpha)=\bigcup\limits_{n=1}^{\infty}{\bigcup\limits_{{\bf i}\in
\Sigma^n}{\left\{{\bf t}\in \RR^k : w^{\bf t}_{\bf i}(A_{\bf
t})\subset \bigcup\limits_{{\bf j}\in \Sigma^n,\ {\bf j}\neq {\bf
i}}{w^{\bf t}_{\bf j}(A_{\bf t})}\right\}}}.
$$
Denote by $\Pi_{\bf k}$, ${\bf k}\in \Sigma^\infty$, the element
${\bf x}$ in $A$ with address ${\bf k}$. Let also $\Pi^{\bf t}_{\bf
k}$, ${\bf k}\in \Sigma^\infty$, be the element in $A_{\bf t}$ with
address ${\bf k}$. For every $n\in \NN$ and ${\bf i}\in \Sigma^n$,
let ${\bf k}({\bf i})\in \Sigma^\infty$ be such sequence that
$$
\Pi_{{\bf i}{\bf k}({\bf i})}\notin \bigcup\limits_{{\bf j}\in
\Sigma^n \atop {\bf j}\neq {\bf i}}{w_{\bf j}(A)}
$$
(such ${\bf k}({\bf i})$ exists since $W$ satisfies the SMP). For
every ${\bf i}\neq {\bf j}\in \Sigma^n$, let
$$
Q_{{\bf i},{\bf j}}=\{{\bf t}\in \RR^k : \Pi^{\bf t}_{{\bf i}{\bf
k}({\bf i})}\in w^{\bf t}_{\bf j}(A_{\bf t})\}.
$$
Then
$$
P({\boldsymbol \alpha})\subset
\bigcup\limits_{n=1}^{\infty}{\bigcup\limits_{{\bf i}\in
\Sigma^n}{\bigcup\limits_{{\bf j}\in \Sigma^n\atop {\bf j}\neq {\bf
i}}{Q_{{\bf i},{\bf j}}}}}.
$$
We now fix a number $n\in \NN$ and indices ${\bf i}\neq {\bf j}\in
\Sigma^n$. For every $m\in \NN$ and ${\bf k}\in \Sigma^m$, denote
$$
Q^{\bf k}_{{\bf i},{\bf j}}=\{{\bf t}\in \RR^k : \Pi^{\bf t}_{{\bf
i}{\bf k}({\bf i})}\in w^{\bf t}_{\bf jk}(A_{\bf t})\}=\{{\bf t}\in
\RR^k : \Pi^{\bf t}_{{\bf i}{\bf k}({\bf i})}=\Pi^{\bf t}_{\bf jkp}\
{\rm for\ some \ {\bf p}\in \Sigma^\infty}\}.
$$
Then
$$
Q_{{\bf i},{\bf j}}=\bigcup\limits_{{\bf k}\in \Sigma^m}{Q^{\bf
k}_{{\bf i},{\bf j}}}.
$$
It is a straightforward argument to verify that for every address
${\bf q}=(q_1,q_2,\ldots)\in \Sigma^\infty$, we have
\begin{align*}
\Pi^{\bf t}_{\bf q}&=\sum\limits_{i=1}^{\infty}{B_{q_1}\cdot \ldots
\cdot B_{q_{i-1}} (I_d-B_{q_i}) (\boldsymbol
\alpha_{q_i}+t_1\boldsymbol \gamma^1_{q_i}+\ldots
+t_k\gamma^k_{q_i})}\\
&=\sum\limits_{i=1}^{\infty}{B_{q_1}\cdot \ldots \cdot B_{q_{i-1}}
(I_d-B_{q_i}) \boldsymbol
\alpha_{q_i}}+t_1\sum\limits_{i=1}^{\infty}{B_{q_1}\cdot \ldots
\cdot B_{q_{i-1}} (I_d-B_{q_i})\boldsymbol \gamma^1_{q_i}}\\
&\quad +\ldots +
t_k\sum\limits_{i=1}^{\infty}{B_{q_1}\cdot \ldots \cdot B_{q_{i-1}}
(I_d-B_{q_i})\boldsymbol \gamma^k_{q_i}}=\Pi_{\bf q}(\boldsymbol
\alpha)+t_1\Pi_{\bf q}(\boldsymbol \gamma_1)+\ldots+t_k\Pi_{\bf
q}(\boldsymbol \gamma_k).
\end{align*}
Then
\begin{align*}
Q^{\bf k}_{{\bf i},{\bf j}}&=\left\{{\bf t}\in \RR^k : \Pi_{{\bf i}{\bf
k}({\bf i})}(\boldsymbol \alpha)+\sum\limits_{i=1}^{k}{t_i\Pi_{{\bf
i}{\bf k}({\bf i})}(\boldsymbol\gamma_i)}=\Pi_{\bf jkp}(\boldsymbol
\alpha)+\sum\limits_{i=1}^{k}{t_i\Pi_{\bf jkp}(\boldsymbol\gamma_i)}
 \mbox{ for some } p\in \Sigma^\infty\right\}\\
&=\left\{{\bf t}\in\RR^k : \sum\limits_{i=1}^{k}{t_i(\Pi_{{\bf i}{\bf
k}({\bf i})}(\boldsymbol\gamma_i)-\Pi_{\bf
jkp}(\boldsymbol\gamma_i))}=\Pi_{\bf jkp}(\boldsymbol
\alpha)-\Pi_{{\bf i}{\bf k}({\bf i})}(\boldsymbol \alpha) \mbox{ 
for some } p\in \Sigma^\infty\right\}.
\end{align*}
Given an address $q\in \Sigma^\infty$, let
$$
B({\bf q})=\[\Pi_{{\bf i}{\bf k}({\bf
i})}(\boldsymbol\gamma_1)-\Pi_{\bf
jq}(\boldsymbol\gamma_1),\ldots,\Pi_{{\bf i}{\bf k}({\bf
i})}(\boldsymbol\gamma_k)-\Pi_{\bf jq}(\boldsymbol\gamma_k)\]
$$
be the $d\times k$ matrix with columns $\Pi_{{\bf i}{\bf k}({\bf
i})}(\boldsymbol\gamma_i)-\Pi_{\bf jq}(\boldsymbol\gamma_i)$,
$i=1,\ldots,k$. Let also ${\bf b}({\bf q})=\Pi_{\bf jq}(\boldsymbol
\alpha)-\Pi_{{\bf i}{\bf k}({\bf i})}(\boldsymbol \alpha)$,
$\sigma=\max_{i=1,\ldots,N}{\|B_i\|}$,
$$
a=\max\{{\rm diam} \ A(\boldsymbol \alpha),{{\rm diam}A(\boldsymbol
\gamma_1)},\ldots,{\rm diam}A(\boldsymbol \gamma_k)\},
$$
where $A({\bf c})$, ${\bf c}=({\bf c}_1,\ldots,{\bf c}_N)$, denotes
the attractor of the system $u_i({\bf x})=B_i({\bf x}-{\bf
c}_i)+{\bf c}_i$, $i=1,\ldots,N$, and for an index ${\bf
j}=(j_1,\ldots,j_n)\in \Sigma^n$, denote $\sigma _{\bf
j}=\|B_{j_1}\|\cdot \ldots \cdot \|B_{j_n}\|$. Then
$$
Q^{\bf k}_{{\bf i},{\bf j}}=\{{\bf t}\in \RR^k : B({\bf
kp})\cdot{\bf t}={\bf b}({\bf kp}) {\rm \ for \ some \ }{\bf p}\in
\Sigma^\infty\}.
$$
We will need the following auxiliary statement.
\begin {lemma}\label {Lemma10}
Let $1\leq k\leq d$ be integers, $\mathcal C$ be a set of $d\times
k$ matrices of rank $k$, which has diameter $\delta$ with respect to
the matrix norm (\ref {norm}), and $\mathcal P$ be a set of vectors
from $\RR^d$, which has diameter $\epsilon$ with respect to the
Euclidean distance. Assume that there exists a finite and positive
number $M>0$ such that for every matrix $B\in \mathcal C$,
$$
\left\|\(B^TB\)^{-1}\right\|\leq M.
$$
Denote also by $L$ and $K$ positive numbers such that $ \|B\|\leq L
$ for every matrix $B\in \mathcal C$, and $ \left|{\bf b}\right|\leq
K $ for every vector ${\bf b}\in \mathcal P$. Let $Q$ be the set of
all vectors ${\bf t}\in \RR^k$, which are solutions to the equation
$$
B{\bf t}={\bf b}
$$
for some matrix $B\in \mathcal C$ and vector ${\bf b}\in \mathcal
P$. Then
\begin {equation}\label {p2'}
{\rm diam}\ Q\leq \epsilon ML+\delta MK + 2\delta M^2L^2K.
\end {equation}
\end {lemma}
{\bf Proof.} Let ${\bf t}_1$ and ${\bf t}_2$ be arbitrary points
from $Q$. There exist matrices $B_1,B_2\in \mathcal C$ and vectors
${\bf b}_1,{\bf b}_2\in \mathcal P$ such that
\begin {equation}\label {p1}
B_i{\bf t}_i={\bf b}_i,\ \ \ i=1,2.
\end {equation}
Since matrices $B_1$ and $B_2$ have rank $k$, each equation (\ref
{p1}) has a unique solution ${\bf
t}_i=\left(B^T_iB_i\right)^{-1}B_i^T{\bf b}_i$, $i=1,2$. Then
\begin{align*}
\left|{\bf t}_1-{\bf
t}_2\right|&=\left|\left(B^T_1B_1\right)^{-1}B_1^T{\bf
b}_1-\left(B^T_2B_2\right)^{-1}B_2^T{\bf b}_2\right|\\
&\leq \left|\left(B^T_1B_1\right)^{-1}B_1^T{\bf
b}_1-\left(B^T_2B_2\right)^{-1}B_1^T{\bf
b}_1\right|+\left|\left(B^T_2B_2\right)^{-1}B_1^T{\bf
b}_1-\left(B^T_2B_2\right)^{-1}B_2^T{\bf
b}_1\right|\\
&\quad + \left|\left(B^T_2B_2\right)^{-1}B_2^T{\bf
b}_1-\left(B^T_2B_2\right)^{-1}B_2^T{\bf b}_2\right|\\
& \leq
\|\left(B^T_1B_1\right)^{-1}-\left(B^T_2B_2\right)^{-1}\|\cdot
\|B_1^T\|\cdot \left|{\bf b}_1\right|\\
&\quad +\|\left(B^T_2B_2\right)^{-1}\|\cdot \|B^T_1-B^T_2\|\cdot \left|{\bf
b}_1\right|+\|\left(B^T_2B_2\right)^{-1}\|\cdot\|B^T_2\|\cdot
\left|{\bf b}_1-{\bf b}_2\right|.
\end{align*}
Due to equality $\|B^T\|=\|B\|$ and definition of numbers $M,L$ and
$K$, we have
$$
\left|{\bf t}_1-{\bf t}_2\right|\leq LK\cdot
\|\left(B^T_1B_1\right)^{-1}-\left(B^T_2B_2\right)^{-1}\|+\delta
MK+\epsilon ML.
$$
Using the following estimate
\begin{align}\label {IB}
\|\left(B^T_1B_1\right)^{-1}-\left(B^T_2B_2\right)^{-1}\|&=
\|\left(B^T_2B_2\right)^{-1}B^T_2B_2\left(B^T_1B_1\right)^{-1}
-\left(B^T_2B_2\right)^{-1}B^T_1B_1\left(B^T_1B_1\right)^{-1}\|\nonumber\\
&=\|\left(B^T_2B_2\right)^{-1}\(B^T_2B_2-B^T_1B_1\)\left(B^T_1B_1\right)^{-1}\|\leq
M^2\|B^T_2B_2-B^T_1B_1\|\nonumber\\
&\leq M^2\(\|B^T_2B_2-B^T_2B_1\|+\|B^T_2B_1-B^T_1B_1\|\)\nonumber\\
&\leq M^2 \(\|B^T_2\|\cdot \|B_2-B_1\|+\|B^T_2-B^T_1\|\cdot
\|B_1\|\)\leq 2\delta M^2L,
\end {align}
for every ${\bf t}_1,{\bf t}_2\in Q$, we obtain
$$
\left|{\bf t}_1-{\bf t}_2\right|\leq 2\delta M^2L^2 K+\delta MK
+\epsilon ML,
$$
and estimate (\ref {p2'}) follows. Lemma \ref {Lemma10} is
proved.\eproof

{\bf Completion of the proof of Proposition \ref {P1}.} We apply
Lemma \ref {Lemma10} with $\mathcal C=\{B({\bf kp}) : {\bf p}\in
\Sigma^\infty\}$ and $\mathcal P=\{{\bf b}({\bf kp}) : {\bf p}\in
\Sigma^\infty\}$. For a matrix $B=[{\bf b}_1,\ldots,{\bf b}_k]$,
denote
$$
\|B\|_{2,\infty}=\max _{i=1,\ldots,k}{\left|{\bf b}_i\right|}.
$$
It is not difficult to see that for any $d\times k$ matrix $B$,
\begin {equation}\label {k'}
\|B\|_{2,\infty}\leq \|B\|\leq \sqrt {k}\|B\|_{2,\infty}.
\end {equation}
Let
$$
M_{{\bf i},{\bf j}}=\sup\limits_{{\bf q}\in \Sigma^\infty}{\|(B({\bf
q})^TB({\bf q}))^{-1}\|}.
$$
Denote $\mathcal Y=\{B({\bf q}) : {\bf q}\in \Sigma^\infty\}$. By
assumption, the columns of matrix $B({\bf q})$ are linearly
independent for every ${\bf q}\in \Sigma^\infty$. In view of the
fact that ${\rm det}B^TB\neq 0$, $B\in \mathcal Y$, and continuity
of ${\rm det}B^TB$ and of the algebraic complement to every element
of $B^TB$, we have that $\|(B^TB)^{-1}\|$ is also continuous with
respect to matrix $B\in \mathcal Y$. Since $\mathcal Y$ is compact
with respect to the matrix norm (\ref {norm}), we obtain that
$M_{{\bf i},{\bf j}}$ is finite.

It is not difficult to see that ${\rm diam}\ \mathcal C\leq a\sqrt k
\sigma_{\bf j}\sigma_{\bf k}$ and ${\rm diam}\ \mathcal P\leq
a\sigma_{\bf j}\sigma_{\bf k}$. Denote
$$
L_{{\bf i},{\bf j}}=\sup\limits_{{\bf q}\in \Sigma^\infty}{\|B({\bf
q})\|},
$$
and let
$$
K_{{\bf i},{\bf j}}=\sup\limits_{{\bf q}\in
\Sigma^\infty}{\left|{\bf b}({\bf q})\right|}.
$$
Then by Lemma \ref {Lemma10},
$$
{\rm diam}\ Q^{\bf k}_{{\bf i},{\bf j}}\leq \sigma_{\bf jk}a M_{{\bf
i},{\bf j}}L_{{\bf i},{\bf j}}+\sigma_{\bf jk}\sqrt k a M_{{\bf
i},{\bf j}}K_{{\bf i},{\bf j}}+2\sigma_{\bf jk}\sqrt k a M_{{\bf
i},{\bf j}}^2L_{{\bf i},{\bf j}}^2K_{{\bf i},{\bf j}}=:\sigma_{\bf
j}\sigma_{\bf k}U_{{\bf i},{\bf j}}.
$$
Denote by $\lambda$ such number that
$$
\sum\limits_{i=1}^{N}{\|B_i\|^\lambda}=1.
$$
Then
$$
\mathcal H_\lambda(Q_{{\bf i},{\bf j}})\leq
\limsup\limits_{m\to\infty}{\sum\limits_{{\bf k}\in \Sigma^m}\({\rm
diam}\ Q^{\bf k}_{{\bf i},{\bf
j}}\)^\lambda}\leq\lim\limits_{m\to\infty}{\sum\limits_{{\bf k}\in
\Sigma^m}{\sigma^\lambda_{{\bf j}}\sigma^\lambda_{\bf
k}U^\lambda_{{\bf i},{\bf j}}}}=\sigma^\lambda_{\bf
j}U^\lambda_{{\bf i},{\bf j}}<\infty.
$$
Since $P(\boldsymbol \alpha)$ is covered by a countable collection
of sets of Hausdorff dimension at most $\lambda$, we have ${\rm dim}
P(\boldsymbol \alpha)\leq \lambda<k$. Hence, the complement of
$P(\boldsymbol\alpha)$ is dense in $\RR^k$ and we can find vector
${\bf t}=(t_1,\ldots,t_k)\in \RR^k$ such that point $\boldsymbol
\alpha_i+\sum\limits_{j=1}^{k}{t_j\gamma^j_i}$ will be arbitrarily
close to $\alpha_i$, $i=1,\ldots,N$, and the system $W_{\bf t}$ will
have the SMP. This implies that $E_d(B_1,\ldots,B_N)$ is dense in
$\RR^{dN}$. Corollary~\ref {C3} implies that $E_d(B_1,\ldots,B_N)$
is a $G_\delta$-set. Proposition~\ref {P1} is proved.\eproof

{\bf Proof of Theorem \ref {Th5}.} Let ${\bf u}\in \RR^d$ be a unit
vector. Since $\sum_{i=1}^{N}{\|B_i\|}<1$, there are numbers
$c_1,\ldots,c_N\in (-1,1)$ such that balls $B[c_i{\bf u},\|B_i\|]$,
$i=1,\ldots,N$, are pairwise disjoint and are contained in $B[{\bf
0},1]$. Let $\boldsymbol\gamma^1_i=c_i(I_d-B_i)^{-1}{\bf u}$, and
$$
w_i({\bf x})=B_i{\bf x}+c_i{\bf u}=B_i({\bf
x}-\boldsymbol\gamma^1_i)+\boldsymbol\gamma^1_i,\ \ i=1,\ldots,N,
$$
and $A=A(w_1,\ldots,w_N)$ be the attractor of the system
$(w_1,\ldots,w_N)$. It is not difficult to see that
$$
w_i(B[{\bf 0},1])\subset B[c_i{\bf u},\|B_i\|]\subset B[{\bf 0},1],\
\ i=1,\ldots,N.
$$
This implies that $A\subset B[{\bf 0},1]$. Indeed, for every element
${\bf x}\in A$, there is a sequence $(i_1,i_2,\ldots)\in
\Sigma^\infty$ such that ${\bf
x}=\lim\limits_{n\to\infty}{w_{i_1\ldots i_n}}({\bf 0})$. Since
$w_{i_1\ldots i_n}({\bf 0})\in B[{\bf 0},1]$ for every $n\in \NN$,
we have ${\bf x}\in B[{\bf 0},1]$. Then we also have
$$
w_i(A)\cap w_j(A)\subset w_i(B[{\bf 0},1])\cap w_j(B[{\bf
0},1])\subset B[c_i{\bf u},\|B_i\|]\cap B[c_j{\bf
u},\|B_j\|]=\emptyset,\ \ i\neq j,
$$
which implies $w_{\bf i}(A)\cap w_{\bf j}(A)=\emptyset$, ${\bf
i},{\bf j}\in \Sigma^n$, ${\bf i}\neq {\bf j}$, $n\in \NN$. Hence,
system of mappings $(w_1,\ldots,w_N)$ has the SMP and we have
$\boldsymbol \gamma_1=(\boldsymbol
\gamma^1_1,\ldots,\boldsymbol\gamma^1_N)\in E_d(B_1,\ldots,B_N)$.
Since $k=1$, the other assumption of Proposition \ref {P1} does not
apply and the density of $E_d(B_1,\ldots,B_N)$ as well as the fact
that it is a $G_\delta$-set follows. Theorem \ref {Th5} is
proved.\eproof

{\bf Proof of Theorem \ref {Th7}.} Assume that
$E_2(B_1,\ldots,B_N)\neq \emptyset$ and let
$\boldsymbol\gamma_1=(\boldsymbol\gamma^1_1,\ldots,\boldsymbol\gamma^1_N)\in
(\RR^2)^N$ be such collection of points that the system of mappings
$$
w_i({\bf x})=B_i({\bf
x}-\boldsymbol\gamma_i^1)+\boldsymbol\gamma^1_i,\ \ i=1,\ldots,N,
$$
satisfies the SMP. Denote
$$
V=\left(
      \begin{array}{cc}
        0\ \   -1 \\
        1\ \  \ \ \  0 \\
      \end{array}
    \right)
$$
($V$ is a rotation matrix) and let
$\boldsymbol\gamma_2=\(V\boldsymbol\gamma^1_1,\ldots,V\boldsymbol\gamma^1_N\)$.
Note that for any non-zero vector ${\bf x}=(x_1,x_2)\in \RR^2$, we
have
$$
{\rm det}\[{\bf x},V{\bf x}\]=\(\left(
      \begin{array}{cc}
        x_1 \ \  -x_2 \\
        x_2\ \ \ \ \ \    x_1 \\
      \end{array}
    \right)
\)=x_1^2+x_2^2\neq 0.$$ 
Since rotation matrices commute, for every
address ${\bf q}=(q_1,q_2,\ldots)\Sigma^\infty$, we obtain,
\begin{align*}
\Pi_{\bf
q}(\boldsymbol\gamma_2)&=\sum\limits_{i=1}^{\infty}{B_{q_1}\cdot\ldots\cdot
B_{q_{i-1}}(I_2-B_{q_i})V\boldsymbol\gamma^1_i}=\sum\limits_{i=1}^{\infty}{B_{q_1}\cdot\ldots\cdot
B_{q_{i-1}}(V-B_{q_i}V)\boldsymbol\gamma^1_i}\\
&=\sum\limits_{i=1}^{\infty}{B_{q_1}\cdot\ldots\cdot
B_{q_{i-1}}(V-VB_{q_i})\boldsymbol\gamma^1_i}=V\sum\limits_{i=1}^{\infty}{B_{q_1}\cdot\ldots\cdot
B_{q_{i-1}}(I_2-B_{q_i})\boldsymbol\gamma^1_i}=V\Pi_{\bf
q}(\boldsymbol\gamma_1).
\end{align*}
Then for every pair of addresses ${\bf i}\neq {\bf j}\in
\Sigma^\infty$ such that $\Pi_{\bf i}(\boldsymbol
\gamma_1)\neq\Pi_{\bf j}(\boldsymbol \gamma_1)$, we have
\begin{align*}
&{\rm det}\ \[\Pi_{\bf i}(\boldsymbol \gamma_1)-\Pi_{\bf
j}(\boldsymbol \gamma_1),\Pi_{\bf i}(\boldsymbol \gamma_2)-\Pi_{\bf
j}(\boldsymbol \gamma_2)\]\\
&\quad ={\rm det}\ \[\Pi_{\bf i}(\boldsymbol \gamma_1)-\Pi_{\bf
j}(\boldsymbol \gamma_1),V(\Pi_{\bf i}(\boldsymbol
\gamma_1)-\Pi_{\bf j}(\boldsymbol \gamma_1))\]\neq 0.
\end{align*}
Then vectors $\Pi_{\bf i}(\boldsymbol \gamma_i)-\Pi_{\bf
j}(\boldsymbol \gamma_i)$, $i=1,2$, are linearly independent and by
Proposition \ref {P1} we obtain that $E_2(B_1,\ldots,B_N)$ is a
dense $G_\delta$-subset of $\RR^{2N}$.\eproof

\end {document}